\documentclass[graybox]{svmult}

\usepackage{type1cm}        
\usepackage{makeidx}         
\usepackage{graphicx}        
\usepackage{multicol}        
\usepackage[bottom]{footmisc}

\usepackage{newtxtext}       %
\usepackage{newtxmath}       

\usepackage[%
   colorlinks=true,
   bookmarks=true,
   unicode,
   plainpages=false,
   ]{hyperref}
\usepackage{url}
   \usepackage{mathtools}
\newcommand{\pot}{U}
\newcommand{\potV}{U}
\newcommand{\Grenze}{g}
\DeclareMathOperator{\spt}{spt}
\newcommand{\function}{\omega}
\newcommand{\multappomega}{\function^{A}}
\newcommand{\appAq}[1]{{#1}^{A}}
\newcommand{\MC}{a}
\newcommand{\param}{\lambda}
\newcommand{\para}{{\delta}}
\newcommand{\diver}{\operatorname{div}}
\newcommand{\dist}{\operatorname{dist}}
\newcommand{\tran}[1]{{#1}_{-\tau}}
\newcommand{\trap}[1]{{#1}_{\tau}}
\newcommand{\difp}[1]{d^+{#1}}
\newcommand{\trapm}[1]{{#1}_{\pm\tau}}
\newcommand{\difpm}[1]{d^\pm{#1}}
\newcommand{\difn}[1]{d^-{#1}}
\newcommand{\intO}{\int\limits_{\Omega}}
\newcommand{\otimess}{\overset{s}{\otimes}}
 \newcommand{\td}{\partial_{\tau}}
\newcommand{\bue}{\bu^{A}}
\renewcommand{\O}{\Omega}

\usepackage[rose,frak,operator]{paper_diening}
\makeindex             

\begin{document}
\title*{Natural second order regularity for systems in the case $1<p\leq 2$ using the
  $A$-approximation}
\titlerunning{Natural second order regularity using the
    $A$-approximation}
\author{Luigi C.\ Berselli  and Michael R\r u\v zi\v cka{}
 }
 \institute{Luigi C.\ Berselli \at Dipartimento di Matematica, Universit{\`a} di Pisa, Via
   F.~Buonarroti 1/c, I-56127 Pisa, ITALY, \email{luigi.carlo.berselli@unipi.it} \and
   Michael R\r u\v zi\v cka{} \at Institute of Applied Mathematics,
   Albert-Ludwigs-University Freiburg, Ernst-Zermelo-Str.~1, D-79104 Freiburg, GERMANY
   \email{rose@mathematik.uni-freiburg.de}}
%
%

\maketitle

\abstract*{In this paper we consider nonlinear problems with an operator depending only on
  the deformation tensor. We consider the class of operators derived from a potential and
  with $(p,\delta)$ structure, for $1<p\leq 2$ and for all $\delta\geq0$. We apply the so
  called $A$-approximation method to approximate the operator by another one with linear
  growth. This allows us to prove the ``natural'' second order regularity (up to the
  boundary) in the case of homogeneous Dirichlet boundary conditions. We focus on the
  steady (elliptic) case, but results are given also in the time-dependent (parabolic)
  case. Results presented are not completely new, but the method we apply was not used
  before in this setting.  }
%
\abstract{In this paper we study nonlinear problems with an operator depending only on
  the symmetric gradient. We consider the class of operators, derived
  from a potential, having $(p,\delta)$-structure for some $1<p\leq 2$ and some $\delta\geq0$. We apply the so
  called $A$-approximation method to approximate the operator by another one with linear
  growth. This allows us to prove ``natural'' second order regularity (up to the
  boundary) in the case of homogeneous Dirichlet boundary conditions. We focus on the
  steady (elliptic) case. However, corresponding results are stated also in the time-dependent (parabolic)
  case. 
}
\section{Introduction}
\label{sec:1}
In this paper we consider the boundary value problem associated to
nonlinear elliptic systems
\begin{align}
  \label{eq:plasticity}
  \begin{aligned}
    -\divo \bfS (\bfD\bfu) &= \bff
    \qquad&&\text{in } \Omega,
    \\
    \bu &= \bfzero &&\text{on } \partial \Omega\,,
  \end{aligned}
\end{align}
where the operator $\bS$ depends on the symmetric gradient
$\bD\bu$ 
and has
$(p,\delta)$-structure (cf.~Definition \ref{def:ass_S}). Here
$\Omega\subset\setR^{3}$ is a sufficiently smooth and bounded domain.
The paradigmatic example for the operator in \eqref{eq:plasticity} is
given by
\begin{equation}
  \label{eq:example}
  \mathbf{S}(\bD\bu):=(\delta+|\bD\bu|)^{p-2}\bD\bu,
  \quad\text{with}\quad \delta\geq0,\ 1<p<\infty\,.
\end{equation}
Thus, problem \eqref{eq:plasticity} is a generalization to systems of the classical
$p$-Laplace problem for scalars $\Delta_{p}u:=\divo(|\nabla
u|^{p-2}\nabla u)$, which corresponds to the case
$\delta=0$. While the existence of weak solutions is a rather standard
result --based on the theory of monotone operators-- the regularity of
solutions is more complicated and has been addressed for the case
$1<p\leq 2$ by Seregin and Shilkin~\cite{SS00} (in the case of a flat boundary)
and by the authors of the present paper in \cite{br-plasticity} (in a
general smooth domain). The proof is obtained by a classical strategy:
the use of difference quotients to estimate partial derivatives
in the tangential directions and ellipticity to recover normal
derivatives. The main difficulties are those of justifying the
calculations to make the argument rigorous. This has been done by
means of: a) smoothing with the addition of an extra Laplace term
$-\epsilon\Delta\bu^\vep$; b) proving for the solution
$\bu^{\epsilon}$ (of the approximate problem) estimates independent of
$\epsilon>0$, to justify the limit as $\epsilon\to0$.

This approach cannot be used in the case
$p>2$, since the calculations --despite being formally very similar--
are not justified. In fact, the
added 
Laplacian term immediately implies estimates in
$L^{2}(\Omega)$ for second order partial derivatives of $\bu^\vep$, but this is
still not enough to give proper meaning to all of the
integrals appearing in the derivation of the various estimates.

To overcome this technical problem --very recently-- we developed
in~\cite{br-multiple-approx} a theory based on the \textit{(multiple) approximation} of
the operator $\bS$, which allows to treat the case $p>2$, for all arbitrarily large $p$.
The theory of the multiple-approximation can also be applied in the case $1<p\leq 2$ (in
fact a single approximation is enough in this case), providing an
alternative proof for the
results from~\cite{SS00,br-plasticity}.

In this paper we consider the case $1<p\leq 2$ and we explain the
modifications and simplifications of the theory with a ``single''
$A$-approximation. Even if the results we prove are not completely
original, we believe it is important to explain them with great
detail. This will be particularly interesting for students or younger
researchers, since the developed method, which is highly flexible, 
can be adapted to several other problems. Even if we
skip some details (which would make the presentation too long) we try
to keep the presentation as much as possible self contained. We refer
with detailed citations
to~\cite{br-reg-shearthin,br-plasticity,br-parabolic, br-multiple-approx} for all
missing technical details.
We present a detailed presentation only in 
the elliptic case. Nevertheless,  the method can be also applied to parabolic
problems with minor modifications, to recover in a different way
results similar to those proved in~\cite{br-parabolic} (see Section \ref{sec:parabolic}).

\medskip

The main goal of this paper is to show how to prove a result of
``natural'' second-order regularity for weak solutions. This
corresponds to proving --under possibly minimal assumptions on the
data-- that weak solutions (and not solutions with additional unproved
properties) satisfy the following inequality
\begin{equation}
  \label{eq:natural}
  \int\limits_ \Omega (\delta+|\bD\bu|)^{p-2}|\nabla\bD\bu|^{2}\,d\bx
  \leq C\,, 
\end{equation}
which can be equivalently rewritten as
$ \nabla \bF(\bD\bu) \in L^{2}(\Omega)$, where 
\begin{equation}
  \label{eq:F}
  \bF(\bD\bu):=(\delta+|\bD\bu|)^{\frac{p-2}{2}}|\bD\bu|\,.
\end{equation}
The regularity coming from inequality \eqref{eq:natural} is called
natural since if one restricts to the periodic-case (and integration
by parts can be done freely without boundary terms) this is formally
obtained by multiplying the system in \eqref{eq:plasticity} by $-\Delta\bu$ and performing straightforward
integration by parts.
%
\begin{remark}
  In the literature the name natural is used to distinguish such
  regularity results from so-called ``optimal'' second-order regularity (for
  which there exists also an intense research activity,
  see~\cite{balci2021pointwise-old,BM2020,CM2019,CM2020}), which proves
  $\nabla \bS\in L^{2}(\Omega)$, i.e., 
  \begin{equation*}
    \int\limits_{\Omega}\big
    |\nabla\big((\delta+|\bD\bu|)^{p-2}\bD\bu\big)\big |^{2}\,d\bx\leq
    C\,.
  \end{equation*}
  The two notions of regularity are rather different in the spirit:
  the optimal regularity is linked with nonlinear versions of the
  singular integral theory, while the natural regularity is based on
  energy methods. It involves quasi-norms (cf.~Barrett and
  Liu~\cite{baliu}), which are, among others, of crucial relevance for
  the numerical analysis of the problem, in particular, to prove
  optimal convergence rates for Finite Element discretizations.
\end{remark}
In Section~2 we will give definitions of the missing notions
and formulate general assumptions on the 
operator $\bS$, covering the example~\eqref{eq:example} in the case
$p\in(1,2]$ and $\delta\in[0, \infty)$. Based on that we consider the
following notion of solution:
\begin{definition}[Regular solution]
  Let the operator $\bS$ in~\eqref{eq:plasticity} have
  $(p,\delta)$-structure for some $p\in(1,\infty)$ and
  $\delta\in[0, \infty)$. 
  We say that $\bu$ is a
  regular solution to \eqref{eq:plasticity} if
  $\bu \in W^{1,p}_0 (\Omega)$ satisfies for all
  $\bw\in W^{1,p}_0(\Omega)$
\begin{equation*}
\begin{aligned}
  \intO{\bS(\bD\bu)}\cdot{\bD\bw}\,d\bx
=
\intO{\bff}\cdot{\bfw}\,d\bx\,,
\end{aligned}
\end{equation*}
and  fulfils
  \begin{align*}
    \begin{split}
      \bF(\bD\bfu)&\in W^{1,2}(\Omega) \,.
    \end{split}
  \end{align*}
\end{definition}
 The main result we will prove with full details is the following:
 \begin{theorem}
  \label{thm:MT}
  Let the operator $\bS$ in~\eqref{eq:plasticity}, derived from a
  potential $\pot $, have $(p,\delta)$-structure for some $p\in(1,2]$,
  and $\delta\in[0, \infty)$. 
  Let $\Omega\subset\setR^3$ be a bounded domain with $C^{2,1}$
  boundary. Assume that $\bff \in L^{p'}(\Omega)$.  Then, the
  system~\eqref{eq:plasticity} has a unique regular solution with
  norms estimated only in terms of the characteristics of $\bS$,
  $\delta$, $\Omega$, and $\|\bff\|_{p'}$.
\end{theorem}
The counterpart in the parabolic case
(cf.~Theorem~\ref{thm:MT-parabolic}) will be presented, without a
detailed proof, in the final section.  Moreover, for all results we
will study only the non-degenerate case $\delta>0$. The degenerate
case can be handled by a limiting argument, provided that the
estimates do not degenerate as $\delta\to0$, exactly as
in~\cite[Sec.~3.2]{br-plasticity}. The proof of such estimates require
some changes with respect to the ones obtained
in~\cite{br-multiple-approx} (for $p>2$) and in~\cite{br-parabolic}
(for $p<2$, but with a different approximation) related to the 
initial condition. Such estimates are available in our setting
(cf.~Proposition~\ref{prop:JMAA2017-1},
Proposition~\ref{prop:main}, Proposition~\ref{prop:main-parabolic1}). In fact, the limiting process
$\delta \to 0$ depends only on the regularity available and is
independent of the method used to prove it, hence there is nothing to
change with respect to the already available proof in
\cite{br-plasticity}.  \vspace{.5cm}

\noindent\textbf{Plan of the paper:} In Section~\ref{sec:section-2} we recall the main
facts about N-functions and the
$A$-approximation. Sections~\ref{sec:3} and \ref{sec:4} are devoted to
the proof of the existence and regularity for the solutions of the
approximated problem. Especially Section~\ref{sec:4} is
crucial for the estimates independent of $A$. Section~\ref{sec:5}
explains the limiting process to come-back from the $A$-approximate
system to the original one. Finally, in Section~\ref{sec:parabolic}
the corresponding results in the parabolic setting are presented.
\section{On the $A$-approximation of an operator and its
  properties} \label{sec:section-2}
In this section we introduce the notation and the crucial properties
of $N$-functions, which will be used to prove the relevant properties of 
$A$-approximated operators. We summarize and recall the main
results already proved with full details
in~\cite{br-multiple-approx,mnr3}, i.e., proofs of all statements in
this section can be found in these references.
\vspace*{-4mm}
\subsection{Notation}
\vspace*{-2mm}
We use $c, C$ to denote generic constants, which may change from line
to line, but are not depending on the crucial quantities. Moreover, we
write $f\sim g$ if and only if there exists constants $c,C>0$ such
that $c\, f \le g\le C\, f$.

We use the customary Lebesgue spaces $(L^p(\Omega), \norm{\,.\,}_p)$,
$p \in [1,\infty]$, and Sobolev spaces
$(W^{k,p}(\Omega), \norm{\,.\,}_{k,p})$, $p \in [1,\infty]$,
$k \in \setN$. 
We do not distinguish between scalar, vector-valued or tensor-valued
function spaces; however, we denote scalar functions by roman letters,
vector-valued functions by small boldfaced letters and tensor-valued
functions by capital boldfaced letters. %
We denote by $\abs{M}$ the $3$-dimensional Lebesgue measure of a
measurable set $M$.
As usual the gradient of a vector field $\bv :\Omega\subset \setR^3 \to \setR^3$ is
denoted as $\nabla \bv = (\partial_{i} v_{j})_{i,j=1,2,3}=(\partial_{i} \bv)_{i=1,2,3}$,
while its symmetric part is denoted as $\bD\bv:= \frac 12 \big (\nabla \bv + \nabla \bv
^\top\big )$.  The derivative of functions defined on tensors, i.e., $\potV:\setR^{3\times
  3}\to \setR$, is denoted as $\partial \potV = (\partial _{i j} \potV)_{i,j=1,2,3}$ where
$\partial _{i j}$ are the partial derivatives with respect to the canonical basis of
$\setR^{3\times 3}$.
%
\vspace*{-4mm}
\subsection{N-functions}
\vspace*{-2mm}
%
A function $\varphi:\setR^{\geq0}\to\setR^{\geq 0}$ is called an {\em
  N-function} 
if $\varphi$ is continuous, convex, strictly positive for $t>0$, and satisfies\footnote{In
  the following we use the convention that $\frac {\varphi'(0)}{0}:=0$.}
  \begin{equation*}
    \lim_{t\to 0^+}\frac{\varphi(t)}{t}=0\,,\qquad\qquad 
    \lim_{t\to \infty}\frac{\varphi(t)}{t}= \infty\,. 
  \end{equation*}
  If $\varphi$ additionally belongs to $C^1(\setR^{\ge 0})\cap
  C^2(\setR^{> 0})$ and satisfies $\varphi''(t)>0$ for all $t>0$, we call
  $\varphi$ a {\em  regular N-function}.  
  In the rest of the paper we restrict ourselves to this case and note
  that for a regular N-function we have $\varphi
  (0)=\varphi'(0)=0$. Moreover, $\varphi'$ is increasing and
  $\lim _{t\to \infty} \varphi'(t)=\infty$. For details we refer
  to~\cite{krasno,Mu,ren-rao,dr-nafsa}. 
%
For a regular N-function $\varphi$ we define the \textit{complementary
  function} $\varphi^*$ via
\begin{equation*}
  \varphi^*(t):=\int\limits_0^t(\varphi')^{-1}(s)\,ds\,.
\end{equation*}
One easily sees that $\varphi^*$ is a regular N-function, too.

  The {\em $\Delta_2$-condition} plays an important role in Orlicz
  spaces. A
 non-decreasing function $\varphi:\setR^{\geq0}\to\setR^{\geq0}$ is said to satisfy the
  {\rm $\Delta_2$-condition} (in short $\varphi\in \Delta_{2}$), if for some constant $K\geq 2$ it holds
  \begin{equation*}
   \varphi(2t)\leq K\varphi(t)\qquad \forall\,t\geq0\,.
\end{equation*}
The $\Delta_2$-constant (the smallest of such $K\geq 2$) of
$\varphi$ is denoted by
$\Delta_2(\varphi)$. 

It has been recently recognized that a fundamental role in regularity
theory of problem similar to \eqref{eq:plasticity} is played by the notion of \textit{balanced N-function}
(cf.~\cite{dr-nafsa,CM2019,br-multiple-approx}). A regular N-function
$\varphi$ is called {\em balanced} if there exist constants $\gamma_1\in (0,1]$ and
$\gamma_2 \ge 1$ such that there holds
  \begin{align*}
    \gamma_1\,\varphi'(t)\le t\,\varphi''(t)\le \gamma_2\,\varphi'(t)
    \qquad  \forall\, t> 0 \,.
  \end{align*}
  The pair $(\gamma_1,\gamma_2)$ is called {\em characteristics} of the balanced
  N-function $\varphi$. The property of being balanced transmits to $\varphi^\ast$, whose
  characteristics is $(\gamma_2^{-1},\gamma_1^{-1})$. %
  Note that for a balanced N-function $\varphi$ we have the equivalences
  \begin{align*}
    \varphi(t)\sim\varphi'(t)\,t \sim  \varphi''(t)\,t^{2}\qquad
    \forall \, t>0
  \end{align*}
  with constants of equivalence depending only on the characteristics
  of $\varphi$. 
%
  In view of this, it is convenient to introduce the particular
  notation $\MC_\varphi: \setR^{\ge 0} \to \setR^{\ge 0}$, defined for regular
  N-functions $\varphi$ via
  \begin{align*}
  \MC_{\varphi}(t):=\frac{\varphi'(t)}{t} \,. 
\end{align*}

Another important tool are {\it shifted N-functions} $\set{\phi_a}_{a \ge 0}$, defined for
$t\geq0$, by
\begin{align*}
  \varphi_a(t):= \int\limits _0^t \varphi_a'(s)\, ds\qquad\text{with }\quad
  \phi'_a(t):=\phi'(a+t)\frac {t}{a+t}\,.
\end{align*}
For an N--function $\phi \in \Delta_2$ there holds for all $\bP,\bQ
  \in \setR^{n\times n}$ and all $t\ge 0$ that $\phi_\abs{\bP}\big
  (\abs{\bP-\bQ}\big )\sim \phi_\abs{\bQ}\big    (\abs{\bP-\bQ}\big )$
  with constants of equivalence depending only on $\Delta_2(\phi')$. The most relevant property for us is a change of
shift.
\begin{lemma}[Change of shift]
  \label{lem:change2}
  Let $\phi$ be an N--function such that $\phi$ and $\phi^*$ satisfy
  the $\Delta_2$--condition.  Then, for all $\delta \in (0,1)$ there
  exists $c_\vep=c_\vep (\Delta_2(\phi'))$ such that for all $\bP,\bQ
  \in \setR^{n\times n}$, and all $t\ge 0$ there holds 
  \begin{align*}
    \phi_{\abs{\bP}}(t) &\le c_\vep\, \phi_{\abs{\bQ}}(t) +
    \vep\, \phi_\abs{\bP}\big    (\abs{\bP-\bQ}\big )\,,\\[1mm]
    \big (\phi_{\abs{\bP}}\big )^*(t) &\le c_\vep\, \big
    (\phi_{\abs{\bQ}}\big )^*(t) +
    \vep\,\phi_\abs{\bP}\big    (\abs{\bP-\bQ}\big )\,.
  \end{align*}
\end{lemma}
\begin{proof}
  These inequalities are proved in~\cite[Lemma 5.15, Lemma 5.18]{dr-nafsa}.
\end{proof}

Finally, we introduce for $p\!\in\!(1,\infty)$ and  $ \para\!\in \!
[0,\infty)$ the function~${\function_{p,\para}:\setR^{\ge 0} \!\to \!
\setR^{\ge 0}}$ via 
\begin{equation*}
 \function_{p,\para}(t):=\int\limits_{0}^{t}(\para+s)^{p-2}s\,ds
 \qquad\forall\,t \ge 0\,.
\end{equation*}
\begin{remark}
  The function $\function_{p,\para}(t)$ is precisely the N-function
  associated with the canonical example for the operator 
 $\mathbf{S}$ in \eqref{eq:example}. If $p$ and $\delta$ are fixed (and
  to avoid confusion with shifted functions) we simply write
  $\function(t):=\function_{p,\para}(t)$.
\end{remark}
Clearly, $\function$ is a regular N-function for all $p\in(1,\infty)$ and $ \para\in
[0,\infty)$. More precisely, for $p\le 2$ we have: 
\begin{lemma}
  \label{lem:function}
  For any $p\in(1,2]$ and for any $\delta\in[0,\infty)$ there holds
  \begin{equation}
    \label{eq:omega}
    \begin{aligned}
     \function(t)&\leq (\function)'(t)\,t\;\leq
      2^{p+1}\function(t)\qquad \forall\,t\geq 0\,,
      \\
      (p-1)  \,    (\function)'(t)&\leq (\function)''(t)\,t\leq
      \,(\function)'(t) \quad \forall\,t >0\,.\hspace*{-5mm}
    \end{aligned}
  \end{equation}
  In particular, 
  the function $\function$ is a balanced N-function with
  characteristics $(p-1,1)$ and $\Delta_2$-constant depending only on
  $p$. Moreover, also $\function^*$ is a balanced N-function with
  characteristics $(1, (p-1)^{-1})$ and $\Delta_2$-constant depending
  only on $p$.
  \end{lemma}
For the shifts of $\function$ and its complementary
function $\function ^*$ there hold for all $a \ge
0$ the equivalences 
$\function_a(t)
\sim (\delta+a+t)^{p-2} t^2$ and $(\function_a)^*(t) \sim \big((\delta+a)^{p-1} +
t\big)^{p'-2} t^2$. 
\vspace*{-1mm}
\subsection{Nonlinear operators with $(p,\delta)$-structure}
\vspace*{-1mm}%
In this section we collect the main properties of nonlinear operators derived from a 
potential and of operators having $(p,\delta)$-structure. 
\begin{definition}[Operator derived from a potential]
  \label{def:potential}
  We say that an operator  \linebreak ${\bS:\,\setR^{3\times3}\to\setR^{3\times
    3}_{\sym}}$ is {\em derived from a potential} 
  $\pot: \setR^{\ge 0} \to \setR^{\ge 0}$, and write $\bS=\partial
  \pot$, if $\bS(\bfzero)=\bfzero$
  and for all $\bP\in\setR^{3\times3}\setminus \set{\bfzero}$ there
  holds
  \begin{align*}
    \bS(\bP)=\partial
    U(|\bP^{\sym}|)=\frac{\pot'(|\bP^{\sym}|)}{\abs{\bP^{\sym}}}\,\bP^{\sym}
    =\MC_{U}(\abs{\bP^{\sym}})\,\bP^{\sym}
  \end{align*}
  for some $\pot\in  C^{1}(\setR^{\geq0})\cap C^{2}(\setR^{>0})$
  satisfying  $\pot(0)=\pot'(0)=0$.
\end{definition}
\begin{definition}[Operator with a  ${\varphi}$-structure]
  \label{def:ass_S}
  Let the operator
  ${\bS\colon \setR^{3 \times 3} \to \setR^{3 \times 3}_{\sym} }$,
  belonging to
  $C^0(\setR^{3 \times 3};\setR^{3 \times 3}_{\sym} )\cap C^1(\setR^{3
    \times 3}\setminus \{\bfzero\}; \setR^{3 \times 3}_{\sym} ) $,
  satisfy ${\bS(\bP) = \bS\big (\bP^{\sym} \big )}$ and
  $\bS(\mathbf 0)=\mathbf 0$. 
  We say that $\bS$ has {\em ${\varphi}$-structure} if there exist 
  a regular N-function $\varphi$ and constants
  $\gamma_3 \in (0,1]$, $\gamma_4 >1$ such that the inequalities
   \begin{equation*}
     \begin{aligned}
       \sum\limits_{i,j,k,l=1}^3 \partial_{kl} S_{ij} (\bP) Q_{ij}Q_{kl}
       &\ge \gamma_3 \, \MC_{\varphi}(\abs{\bfP^{\sym}})\, |\bP^{\sym}
       |^2\,,
       \\
       \big |\partial_{kl} S_{ij}({\bP})\big |
       &\le \gamma_4  \, \MC_{\varphi}(\abs{\bfP^{\sym}})\,,
     \end{aligned}
   \end{equation*}
   are satisfied for all $\bP,\bQ \in \setR^{3\times 3} $ with
   $\bP^{\sym} \neq \bfzero$ and all $i,j,k,l=1,2, 3$.  The constants
   $\gamma_3$, $\gamma_4$, and $\Delta_2(\varphi)$ are called the {\em
     characteristics} of $\bfS$ and will be denoted by
   $(\gamma_3,\gamma_4, \Delta_2(\varphi))$.

   If $\varphi= \function_{p,\delta} $ with $p \in (1, \infty)$ and $\para\in [0,\infty)$
   we say that $\bS$ has {\em $(p,\para)$-structure} and call $(\gamma_3,\gamma_4,p)$ its
   characteristics.
\end{definition}
Closely related to an operator with $\varphi$-structure is 
the function
${\bF_\varphi\colon\setR^{3 \times 3} \to \setR^{3 \times 3}_{\sym}}$ defined via
\begin{equation}\label{def:F}
  \bF_\varphi(\bP):=\frac{\sqrt{\varphi'(\abs{\bfP^{\sym}})\abs{\bP^{\sym}}}}{\abs{\bfP^{\sym}}}\,\bP^{\sym}=\sqrt{\MC_{\varphi}(\abs{\bfP^{\sym}})}\,\bP^{\sym}
  \,,
\end{equation}
where the first representation holds only for $\bP^{\sym} \neq
\bfzero$. 
In the special case of an operator $\bS$ with $(p,\para)$-structure we have (recall that
$\function=\function_{p,\para}$)
\begin{equation*}
  \bF(\bP):=\bF_{\function}(\bP)=\sqrt{\MC_{\function(\abs{\bfP^{\sym}})}}\,\bP^{\sym}
  = \big (\para +\abs{\bP^{\sym}}\big )^{\frac {p-2}2}\bP^{\sym}\,,
\end{equation*}
which is consistent with the notation used in the previous literature,
cf.~\eqref{eq:F}.

If $\varphi$ is a balanced N-function with characteristics
$(\gamma_1,\gamma_2)$, then $\bS=\partial\varphi$ is an operator with
$\varphi$-structure and with characteristics depending only on
$\gamma_1$ and $\gamma_2$.  The following 
result will be crucial for our investigations (cf.~\cite[Section~6]{dr-nafsa}).

\begin{proposition}\label{prop:hammer-phi}
  Let
  $\varphi$ be a balanced N-function with characteristics
  $(\gamma{_1},\gamma{_2})$. Let $\bS$ have
  $\varphi$-structure with characteristics $(\gamma_3,\gamma{_4},
  \Delta{_2(\varphi)})$ and let
  $\bF_\varphi$ be defined in \eqref{def:F}. Then, we have for all
  $\bP,\bQ \in \setR^{3\times 3} $ that
  \begin{align*}
    \big(\bS(\bP)-\bS(\bQ)\big)\cdot(\bP-\bQ)
    &\sim \MC_{\varphi}(\abs{\bfP^{\sym}} +
      \abs{\bfP^{\sym}-\bQ^{\sym}})\,|\bP^{\sym}-\bQ^{\sym}|^{2}
    \\
    &      \sim |\bF_\varphi(\bP)-\bF_\varphi(\bQ)|^{2}\,,
    \\
    |   \bS(\bP)-\bS(\bQ)|&\sim \MC_{\varphi}(\abs{\bfP^{\sym}}+\abs{\bfP^{\sym}-\bfQ^{\sym}} 
                            )\, |\bP^{\sym}-\bQ^{\sym}|\,,
\end{align*}
where the constants of equivalence depend only on 
$\gamma_{1}, \gamma{_2}, \gamma{_3}$, and $\gamma_{4}$.
\end{proposition}
In addition, the following result will be used to handle operators derived from a
potential.
\begin{proposition}
\label{prop:potential-equivalence}
Let the operator $\bS=\partial\potV$, derived from the potential
$\potV$, have \mbox{$\varphi$-structure,} with characteristics
$(\gamma_3,\gamma_4,\Delta_2(\varphi))$. If $\varphi$ is a balanced
N-function with characteristics $(\gamma_1,\gamma_2)$, then $\potV$ is
a balanced N-function satisfying for all $t>0$
\begin{align*}
     \frac{ \gamma_3 }{\gamma_2}\varphi''(t)\le \potV''(t)\le
     \frac{ \gamma_4 }{\gamma_1}\varphi''(t)\,.
  \end{align*}
  The characteristics of $\potV$ is equal to $\big(\frac{\gamma_3}{\gamma_4}
    \,\frac{ \gamma_1^2}{\gamma_2}, \frac{\gamma_4}{\gamma_3}
    \,\frac{ \gamma_2^2}{\gamma_1}\big )$.
\end{proposition}
The significance of this proposition is that a general operator $\bS$ derived from a
potential $\pot$ with $(p,\para)$-structure, can be simply handled as the explicit
example~\eqref{eq:example}. 
\subsection{Approximation of a nonlinear operator}\label{sec:approx}
We now define the $A$-approximation of a function and of an operator,
and prove the relevant properties needed in the sequel. This
approximation was introduced in~\cite{mnr3} for $p>2$, and generalized
in the recent paper~\cite{br-multiple-approx} to a so called
$(A,q)$-approximation, for some $q\geq2$, which allows for a unified
approach for all $p \in (1,\infty)$.  The purpose of the
$A$-approximation of a function is to have
quadratic behaviour near infinity (cf.~\cite[Lemma~2.22]{mnr3}) and
consequently, one can take advantage of the standard Hilbertian theory.
\begin{definition}[$A$-approximation of a
  scalar real function]\label{def:A-approx}  Given a function $\potV
  \in C^{1}(\setR^{\geq0})\cap C^{2}(\setR^{>0})$ satisfying
  $\potV(0)=\potV'(0)=0$, we define for $A\ge 1$ the \smash{\em
    $A$-approximation}
  $\appAq{\potV} \in C^{1}(\setR^{\geq0})\cap C^{2}(\setR^{>0}) $ via
  \begin{equation*}
    \appAq{\potV}(t):=\left\{
      \begin{aligned}
        &\potV(t)\qquad &t\leq A\,,
        \\
        &\alpha_{2}\,t^{2}+\alpha_{1}\,t +\alpha_{0}\qquad &t> A\,.
      \end{aligned}
    \right.
  \end{equation*}
  To ensure continuity up to second order derivatives, the constants
  $\alpha_{i}=\alpha_{i}(\potV)$, $i=0,1,2$, are given via 
  \begin{equation*}
    \begin{aligned}
      \alpha_{2}&=\frac{1}{2}\,\potV''(A)\,,
      \\
      \alpha_{1}&=\potV'(A)-\,\potV''(A)\,A\,,
      \\
      \alpha_{0}&=\potV(A)-\potV'(A)\,A+\frac{1}{2}\potV''(A)\,A^{2}\,.
    \end{aligned}
  \end{equation*}
\end{definition}
\begin{remark}\label{rem:VA}
  If $\varphi$ is a regular N-function, the definition of
  $\varphi^{A}$, together with the properties of $\varphi$, imply that
  there exists a constant $c(A,\varphi)$ such that for all $t\ge 0$
  there holds
  \begin{align*}
      a_{\varphi^{A}}(t) = \frac{(\varphi^{A})'(t)}{t} \le c(A,\varphi)\,.
  \end{align*}
  More precise (explicit) upper and lower bounds are given
  in~\eqref{eq:wAl1} if $\varphi=\function$.
\end{remark}
We have the following relevant result
(cf.~\cite[Lemma~2.42]{br-multiple-approx}) linking balanced functions
with their $A$-approximations.
\begin{lemma}
  \label{lem:eq12}
  Let $\varphi$  be a balanced N-function with characteristics
  $(\gamma_1,\gamma_2)$. Then, for
  all $A\ge 1$,  it holds that $\varphi^{A}$
  is also balanced with characteristics
  $\big(\gamma_1,\gamma_2\big )$. 
\end{lemma}
Concerning the homogeneity of the function
$\function$ for $1<p\leq2$ (similar results  could be  deduced 
also in the case $p>2$)
we have the following result:
\begin{lemma}
  Let $1<p\leq2$ and $\delta\in[0,\infty)$.  The functions
  $\function(t)$ 
  and $ \function^{A}(t)$, for any $A\geq1$, are balanced functions with characteristics
  $(p-1,1)$. Moreover, it holds for all $\lambda,\,t\geq0$ that 
  \begin{equation}
    \label{eq:lambda}
    \function(\lambda\, t)\leq
    \max\{\lambda,    \lambda^{2}\}\,\function(t)\qquad\text{and}\qquad
    \function^{A}(\lambda\, t)\leq \max\{\lambda,
    \lambda^{2}\}\,\function^{A}(t)\,. 
  \end{equation}
\end{lemma}
\begin{proof}
  The assertion on the characteristics for both functions follows directly from Lemma~\ref{lem:function}
  and Lemma~\ref{lem:eq12} (which shows that they are unchanged by the~$A$-~approximation). The estimates in \eqref{eq:lambda} are proved by observing that if $\phi$ is a regular
  N-function with characteristics $(\gamma_{1},\gamma_{2})$, then it follows for all $t>0$
  that
\begin{equation*}
\frac{d}{dt}\log(\varphi'(t))=  \frac{\varphi''(t)}{\varphi'(t)}\leq
\gamma_2\frac{1}{t} \,, 
\end{equation*}
which implies, by integration with respect to $t$ over $(s,\lambda s)$,
with $\lambda>1$ and $s>0$, and using the
exponential function, that
\begin{align*}
  \frac{\varphi'(\lambda\, s)}{\varphi'(s)} \le \lambda^{\gamma_2}\,.
\end{align*}
A further integration with respect to $s$ over $(0,t)$, $t>0$,  proves
\begin{equation*}
  \varphi(\lambda\, t)\leq\lambda^{\gamma_2+1}\varphi(t)\qquad
  \forall\,t>0\,. 
\end{equation*}
The case  $t=0$ is trivial and in the case $0\leq\lambda\leq1$ and
$t\geq0$ the proof ends by observing that
$\phi\big((1-\lambda)\,0+\lambda\, t\big)\leq\lambda\phi(t)$, by the 
convexity of $\phi$.
\end{proof}
Next, we define the $A$-approximation of an operator derived
from a potential.
\begin{definition}[$A$-approximation of an operator derived from a
  potential]\label{def:SA}
  Let the operator $\bS=\partial\pot$ be derived from the potential
  $\pot$.  Then, we define for given $A\ge 1$  the {\em
    $A$-approximation} $\appAq{\bS}:=\partial \appAq{\pot}$ as the operator derived from the potential $\appAq{\pot}$,
  i.e., $\appAq\bS$ satisfies $\appAq\bS(\bfzero)=\bfzero$ and for all
  $\bP\in \setR^{3\times 3}\setminus\set{\bfzero}$ there holds
  \begin{equation*}
    \appAq\bS(\bP):=\partial\appAq{\pot}(|\bP^{\sym}|)
    =\frac{(\appAq\pot)'(|\bP^{\sym}|)}{\abs{{\bP^{\sym}}}}\,\bP^{\sym}
    = \MC_{\appAq{\pot}}(\abs{\bP^{\sym}})\,\bP^{\sym}
    \,.
  \end{equation*}
  \end{definition}
The properties of the operator $\bS$ in
Proposition~\ref{prop:hammer-phi} are inherited by the operator
$\bS^A$. More precisely, we have (cf.~\cite[Prop.~2.47]{br-multiple-approx}):
\begin{proposition}
\label{prop:SA-ham}
Let
$\varphi$ be a balanced N-function with characteristics
$(\gamma_1,\gamma_2)$. Let the operator
$\bS=\partial\potV$, derived from the potential
$\potV$, have
$\varphi$-structure with characteristics
$(\gamma_3,\gamma_4,\Delta_2(\varphi))$.  For $A\ge
1$ let $\appAq{\varphi} $ and $\appAq{\bS} $ be the
$A$-approximation of $\varphi$ and
$\bS$, respectively. %
Then, we have for all
$\bP,\bQ \in \setR^{3\times 3} $ that
\begin{align*}
    (\appAq{\bS}(\bP)-\appAq{\bS}(\bQ))\cdot(\bP-\bQ)
    &\sim \MC_{\phi^{A}}(\abs{\bfP^{\sym}} +
      \abs{\bfP^{\sym}-\bfQ^{\sym}})\,|\bP^{\sym}-\bQ^{\sym}|^{2}  \,,
    \\
    &      \sim |{\bF}_{\phi^{A}}(\bP)-{\bF}_{\phi^{A}}(\bQ)|^{2}, 
    \\
    |   \appAq{\bS}(\bP)-\appAq{\bS}(\bQ)|&\sim \MC_{\phi^{A}}(\abs{\bfP^{\sym}} +
        \abs{\bfP^{\sym}-\bfQ^{\sym}})\,|\bP^{\sym}-\bQ^{\sym}| 
\end{align*}
with constants of equivalence depending only on
$\gamma_1, \gamma_2, \gamma_{3}$, and $\gamma_{4}$.
\end{proposition}
\begin{remark}
  For the  limiting process, it is of fundamental relevance that
  in Proposition~\ref{prop:SA-ham} the constants do not depend on
  $A\geq1$. The details of the proof can be found in \mbox{\cite[Sec.~2]{br-multiple-approx}.} 
\end{remark}

Equivalent expressions for $\nabla\bF(\bD\bu)$ (based on
Proposition \ref{prop:SA-ham}) play a crucial role in the proof of
regularity of weak solutions. To this end, we define, for a sufficiently smooth operator
$\bS:\,\setR^{3\times3}\to\setR^{3\times 3}_{\sym}$, the functions
$\mathbb{P}^{A}_i\colon \setR^{3\times 3} \to \setR$,
$i=1,2,3$, via 
\begin{equation*}
  \begin{aligned}
  \mathbb{P}^{A}_i (\bP)&:=\partial_i\bS^{A}(\bP)\cdot\partial_i\bP =
  \sum_{j,k,l,m=1}^3\partial_{jk}S^{A}_{lm }(\bP) \,\partial_i P_{jk}\,\partial_i P_{lm }\,,
  \end{aligned}
\end{equation*}
and emphasize  that there is no summation over the index $i$. 

The above  mentioned properties of balanced functions allow us to deal with the problem
associated with $\bS=\partial \potV$ with $(p,\para)$-structure using the quantities
related with $\function$ and not with $\potV$ itself, greatly
simplifying both the presentation and the estimates,
cf. Proposition~\ref{prop:potential-equivalence}. For this reason, we also introduce the
following notation, consistent with~\eqref{eq:F}:
\begin{equation*}
  \bF^A(\bP):=\bF_{\function^A}(\bP)\qquad\text{and}\qquad
  \MC^{A}(t):=\MC_{\function^{A}}(t)\,. 
  \end{equation*}

Using this notation we have the following result
(cf.~\cite[Prop.~2.49]{br-multiple-approx}, \cite[Prop.~2.4]{br-plasticity}):
\begin{proposition}\label{prop:pFA}
  Let the operator $\bS= \partial\pot$, derived from the potential
  $\potV$, have $(p,\delta)$-structure for some $p \in (1, \infty)$
  and $\para\in [0,\infty)$, with characteristics
  $(\gamma_3,\gamma_4,p)$. 
   If for a vector field $\bv \colon \Omega \subset \setR^3\to \setR^3$ there
  holds ${\bF^{A}(\bD\bv) \in W^{1,2}(\Omega)}$, then we have for
  $i=1,2,3$ and a.e.~in $\Omega$ the following equivalences
  \begin{align*}
    \begin{aligned}
      |\partial_i \bF^{A}(\bD\bv)|^{2} &\sim\MC^{A}(|\bD\bv|)\,
      |\partial _i\bD\bv|^{2}
      \\
      &\sim \mathbb P^{A}_i(\bD\bv)\,,
      \\
      |\partial_i \bS^{A}(\bD\bv)|^{2} &\sim\MC^{A}(|\bD\bv|)\, \mathbb
      P^{A} _i(\bD\bv)\,,
    \end{aligned}
  \end{align*}
where the constants of equivalence depend only on
$
\gamma_3,\gamma_4$, and $p$.
\end{proposition}

In view of this proposition it is important to have upper and lower bounds for
$\MC^A$ in order to control various quantities
related to $\bF^{A}$, in terms of $\bF$. Crucial in this respect is
the following result (cf.~\cite[Lem.~2.69]{br-multiple-approx}): 
\begin{lemma}\label{lem:UA}
  For $p\in (1,2],\delta> 0$, and $A\ge 1$ the function  $\MC^{A}(t)$ is
  non-increasing and for all $t\ge 0$ there holds
  \begin{align}
    \label{eq:wAl1}
    \begin{gathered}
      (p-1)\, \MC (t) \le \MC^A(t) \le
      \delta^{p-2}\,,
      \\
      (p-1)\, (\delta+A)^{p-2} \le \MC^A(t) \,.
    \end{gathered}
  \end{align}
\end{lemma}
\begin{proof}
  The statement is clear for $t \le A$ using
  $\MC^A(t)=\MC(t)=(\para +t)^{p-2} $, $0\le \para$, $t \le A$, and $p \le 2$. For $t\ge A$ we
  have
  $\MC^A(t)=\function ''(A) + \frac{\function'(A)-\function''(A)\,
    A}{t}$. Thus, we get that $\MC^A(A)=(\delta+A)^{p-2}$,
  $\lim_{t\to
    \infty}\MC^A(t)=(\delta+A)^{p-3}\,\big(\delta+(p-1)A\big)$ and
  $(\MC^A)'(t)=- \frac{\function'(A)-\function''(A)\, A}{t^2}\le 0$ in
  view of~\eqref{eq:omega}, and $p\le 2$.  This yields
  \begin{align*}
    (\para +A)^{p-2}\ge \MC^A(t) \ge (\para
    +A)^{p-3}((p-1)A+ \para) \ge (p-1)\, (\para +A)^{p-2}\,,
  \end{align*}
  which implies the assertions using
  $\delta^{p-2} \ge (\para +A)^{p-2}$ and
  $(\para +A)^{p-2}\ge (\para +t)^{p-2}$ in view of $t\ge A$, and
  $p\le 2$.
\end{proof}
In the sequel we will use frequently the following consequences (cf.~\cite[Cor.~2.71]{br-multiple-approx}):
\begin{corollary}\label{cor:UA}
  Let the operator $\bS$, derived from the potential $\pot$, have
  $(p,\para)$-structure for some $p\in (1,2]$ and $\delta> 0$, with
  characteristics $(\gamma_3,\gamma_4,p)$. %
  Then, there
  holds for all $t\ge 0$ that
  \begin{align*}
    \begin{gathered}
\frac{(p-1)}{2}\, (\delta+A)^{p-2}t^{2} \le \function^A(t) \,,
      \\
      (p-1)\, \function (t) \le \function^A(t) \le
      \frac {\para^{p-2}}2\, t^2\,,
      \\
    (\function^A)^*(t)\le (p-1) \,(\Delta_2(\function^*))^M \,\function^*(t)\,,
    \end{gathered}
  \end{align*}
  where $M\in \setN_0$ is chosen such that $(p-1)^{-1}\le 2^M$. Moreover, 
  for all $\bP \in \setR^{3\times 3}$ there holds
  \begin{align*}
    \begin{aligned}
      \abs{\bF^A(\bP)}^2&\sim \function^A(\abs{\bP^{\sym}})\,,
      \\
    c\, \abs{\bF(\bP)}^2 &\le  \abs{\bF^A(\bP)}^2\,, 
    \\
    \abs{\bS^A(\bP)} &\le c\, \para^{p-2} \abs{\bP^{\sym}} \,
    \end{aligned}
  \end{align*}
  with constants $c$ depending only on $\gamma_3, \gamma_4$, and $p$. 
\end{corollary}
 \begin{corollary}\label{cor:A}
   Under the assumptions of Proposition~\ref{prop:pFA} there exists      $c(p,\gamma_{i})>0$ such that
   \begin{align*}
     c(p,\gamma_{i})    |\nabla \bF(\bD\bv)|^{2}\leq      |\nabla
     \bF^{A}(\bD\bv)|^{2}. 
   \end{align*}
 \end{corollary}
 \begin{proof}
   This follows immediately from Proposition \ref{prop:pFA}, Lemma 
   \ref{lem:UA} and \cite[Prop.~2.4]{br-plasticity}.
 \end{proof}
 \section{On the existence and uniqueness of regular solutions for the   approximate problem}
\label{sec:3}
In this section we introduce the approximate problem and prove
existence, uniqueness, and regularity of its solutions. In fact, to prove
Theorem~\ref{thm:MT} we use an approximate problem, obtained by
replacing the operator $\bS=\partial \pot$ with $(p,\delta)$-structure
by $\bS^{A}= \partial \pot^A$ %
which has $(2,\delta)$ structure, i.e., we study 
\begin{align}
  \label{eq:plasticity-app}
  \begin{aligned}
-\divo \bfS^{A} (\bfD\bue) &= \bff
\qquad&&\text{in }
    \Omega\,,
    \\
    \bue &= \bfzero &&\text{on } \partial \Omega\,.
  \end{aligned}
\end{align}
This system can be treated by standard techniques typical for linear
equations. This procedure yields various estimates
independent of $A\geq1$ for the solution $\bue$, which will enable us
to pass to the limit $A\to\infty$, and to show that the limit
$\bu=\lim_{A\to\infty}\bue$ will be a regular solution of the original
problem \eqref{eq:plasticity}. 

The first standard result concerns the existence and uniqueness of
weak solutions for \eqref{eq:plasticity-app}.
\begin{proposition}
\label{thm:existence_perturbation}
Let the operator $\bS=\partial U$, derived from the potential $\pot$,
have $(p,\delta)$-structure for some $p\in(1,2]$ and
$\delta\in(0,\infty)$. Assume that $\bff \in L^{p'}( \Omega)$. Let
$\bS^{A}$ be 
as in Definition~\ref{def:SA}. Then, the approximate
problem~\eqref{eq:plasticity-app} possesses a unique weak solution,
i.e., 
$\bue\in W^{1,2}_{0}(\Omega)$ with $\bF^A(\bD\bue) \in L^{2}(\Omega)$
satisfies for all 
 $\bw\in W^{1,2}_0(\Omega)$
\begin{equation}
  \label{eq:weak-app}
\begin{aligned}
  \intO{\bS^{A}(\bD\bue)}\cdot {\bD\bfw}\,d\bx 
=
\intO{\bff}\cdot{\bfw}\,  d\bx\,.
\end{aligned}
\end{equation}
This solution  satisfies 
the estimate
\begin{gather}
  \label{eq:main-apriori-estimate2}
  \begin{aligned}
   &\|\bF^{A}(\bD\bue)\|_{2}^{2} +(p-1)(\para+A)^{p-2}\norm{\bD\bue}_2^{2}
    \\
    &\quad
    + (p-1) \big(\|\bF(\bD\bue)\|_{2}^{2} +\norm{\bD \bue}_{p}^{p}\big)
    \leq C\int\limits_{\Omega}\function^{*}(|\bff|)\,d\bx
  \end{aligned}
 \end{gather}
with $C$ depending only on the characteristics of $\bS$ and
$\Omega$.
\end{proposition}
\begin{remark}
  The energy-type estimate \eqref{eq:main-apriori-estimate2}, which is obtained by testing
  with $\bue$, 
  implies that: (i)  $\bue\in W^{1,2}_{0}(\Omega)$
  with norms depending on $A$; \smash{(ii) $\bue \in W^{1,p}_{0}(\Omega)$} with norms
  bounded uniformly with respect to $A\geq1$.
\end{remark}
\begin{proof}[of Proposition~\ref{thm:existence_perturbation}]
  The proof is based on a classical Faedo-Galerkin approximation of
  \eqref{eq:plasticity-app}. The existence of Galerkin solutions $\bue_{k}$, for $k \in \setN$,
  follows by a standard argument based on Brouwer fixed point
  theorem.
  Passing to the limit as
  $k\to \infty$ (for $A$ fixed) can be done within the standard theory
  of monotone operators (Minty-Browder theory). Since this is a
  fully standard argument, we just derive the a priori estimates
  necessary for this procedure.

  By using $\bue _k$ as test function in the Galerkin approximation
  for $\bue_{k} $ we get\\[-1mm]
\begin{equation*}
  \begin{aligned}
c\, \|\bF^{A}(\bD\bue_{k})\|_2^{2} 
    & \leq c_{\epsilon}\int\limits_{\Omega}(\multappomega)^{*}(|\bff|)
    \,d\bx+\epsilon\int\limits_{\Omega}\multappomega(|\bue_{k}|)\,d\bx
    \\
    &\leq c_{\epsilon}\int\limits_{\Omega}(\multappomega)^{*}(|\bff|)
    \,d\bx+\epsilon\,
    C\int\limits_{\Omega}\multappomega(|\bD\bue_{k}|)\,d\bx\,,
  \end{aligned}
\end{equation*}
where we used in the first line Proposition~\ref{prop:SA-ham} with
$\bQ=\bfzero$ together with Young inequality,  and in the second line\\[-1mm]
\begin{equation*}
  \int\limits_{\Omega} \multappomega(|\bue_{k}|)\,d\bx\leq
  C_{P}\int\limits_{\Omega}\multappomega(|\nabla\bue_{k}|)\,d\bx\leq
  C_{P}C_{K}\int\limits_{\Omega}\multappomega(|\bD\bue_{k}|)\,d\bx\,, 
\end{equation*}
which follows from modular versions of Poincar\'e and Korn
inequalities in Orlicz spaces, see~\cite{bdr-phi-stokes,
  br-multiple-approx,talenti}. Moreover, we absorb the last term on  
the right-hand side of the previous estimate using
$\int_{\Omega}\multappomega(|\bD\bue_{k}|)\,d\bx \sim
\|\bF^{A}(\bD\bue_{k})\|_2^{2}$ in view of 
Corollary~\ref{cor:UA}. Note that all constants are independent of
$A\geq1$, and depend only on the characteristics of
$\bS$ and on $\Omega$. Moreover, from Corollary \ref{cor:UA},  it also
follows that
  \begin{equation}\label{eq:f}
    \int\limits_{\Omega}(\multappomega)^{*}(|\bff|)\,d\bx\,
    \leq c(p)    \int\limits_{\Omega}\function^{*}(|\bff|)\,d\bx\leq
    C(p)\Big(\delta^{p}+\int\limits_{\Omega}|\bff|^{p'}\,d\bx\Big)\,, 
  \end{equation}
  where the last estimate shows that the right-hand side in
  \eqref{eq:main-apriori-estimate2} is finite. Hence, after the
  limiting procedure $k \to \infty$ we arrive at
\begin{equation}\label{eq:apri}
  \norm{\bF^{A}(\bD\bue)}_2^{2}
  \leq C\,  \int\limits_{\Omega}\function^{*}(|\bff|)\,d\bx
\end{equation}
for some $C$ independent of $A$ and $\delta$. Uniqueness  follows from the
strict monotonicity of $\bS^{A}$ (cf.~Proposition \ref{prop:SA-ham}). By using the estimates
in Corollary~\ref{cor:UA} and the definition of $\bF^{A}$ we derive
from \eqref{eq:apri} the various terms in the estimate
\eqref{eq:main-apriori-estimate2}, which ends the
proof. 
\end{proof}
\subsection{Description and properties of the boundary}
\label{sec:bdr} 
We assume that the boundary $\partial\Omega$ is of class $C^{2,1}$, that
is for each point $P\in\partial\Omega$ there are local coordinates such
that in these coordinates we have $P=0$ and $\partial\Omega$ is locally
described by a $C^{2,1}$-function, i.e.,~there exist
$R_P,\,R'_P \in (0,\infty),\,r_P\in (0,1)$ and a $C^{2,1}$-function
$\Grenze_{P}:B_{R_P}^{2}(0)\to B_{R'_P}^1(0)$ such that
\begin{itemize}
\item[\rm (b1)]\quad  $\bx\in \partial\Omega\cap (B_{R_P}^{2}(0)\times
  B_{R'_P}^1(0))\ \Longleftrightarrow \ x_3=\Grenze_{P}(x_1,x_2)$\,,
\item   [\rm (b2)]\quad $\Omega_{P}:=\{(x',x_{3})\fdg x'=(x_1,x_2)
 \in  B_{R_P}^{2}(0),\ \Grenze_{P}(x')<x_3<\Grenze_{P}(x')+R'_P\}\subset \Omega$, 
\item [\rm (b3) ]\quad
  $\nabla \Grenze_{P}(0)=\bfzero,\text{ and
  }\forall\,x'=(x_1,x_2)^\top \in B_{R_P}^{2}(0)\quad |\nabla
  \Grenze_{P}(x')|<r_P$\,,
\end{itemize}
where $B_{r}^k(0)$ denotes the $k$-dimensional open ball with center
$0$ and radius ${r>0}$. We also define, for $0<\lambda<1$, the open sets $\lambda\,
\Omega_P\subset \Omega_P$ as 
\begin{equation*}
  \lambda\, \Omega_P:=\{(x',x_{3})\fdg x'=(x_1,x_2)^\top
 \in
  B_{\lambda R_P}^{2}(0),\ \Grenze_{P}(x')<x_3<\Grenze_{P}(x')+\lambda R_P'\}\,.
\end{equation*}
To localize near  $\partial\Omega\cap \partial\Omega_P$, for $P\in\partial\Omega$, we fix smooth
functions $\xi_{P}:\setR^{3}\to\setR$ such that 
\begin{itemize}
\item[$\rm (\ell 1)$]\quad  $\chi_{\frac{1}{2}\Omega_P}(\bx)\leq\xi_P(\bx)\leq
  \chi_{\frac{3}{4}\Omega_P}(\bx)$\,,
\end{itemize}
where $\chi_{A}(\bx)$ is the indicator function of the measurable set
$A$. 
For the remaining interior estimate we localize by a smooth function
${0\leq\xi_{0}\leq 1}$ with $\spt \xi_{0}\subset\Omega_{0}$, where
$\Omega_{0}\subset \Omega$ is an appropriate open set such that
$\dist(\partial\Omega_{0},\,\partial\Omega)>0$.
Since the boundary $\partial\Omega $ is compact, we can use an appropriate
finite sub-covering which, together with the interior estimate, yields
the global estimate.

Let us introduce the tangential derivatives near the boundary. To
simplify the notation we fix $P\in \partial\Omega$, $h\in (0,\frac{R_P}{16})$,
and simply write $\xi:=\xi_P$, $\Grenze:=\Grenze_{P}$. We use the standard notation
$\bx =(x',x_3)^\top$ and denote by $\be^i,i=1,2,3$ the canonical
orthonormal basis in $\setR^3$. In the following lower-case Greek
letters take values $1,\, 2$. For a function $f$ with $\spt
f\subset\spt\xi$ we define for $\alpha=1,2$ tangential translations:
\begin{equation*}
\begin{aligned}
  \trap f(x',x_3) = f_{\tau _\alpha}(x',x_3)&:=f\big (x' +
  h\,\be^\alpha,x_3+\Grenze(x'+h\,\be^\alpha)-\Grenze(x')\big )\,,
\end{aligned}
\end{equation*}
tangential differences $\Delta^+ f:=\trap f-f$, and tangential
difference quotients
%
  $\difp f:= h^{-1}\Delta^+ f$.  
%
For
simplicity we denote $\nabla \Grenze:=(\partial_1 \Grenze,\partial_{2}\Grenze, 0)^\top$
and use the operations $\trap {(\cdot)}$, $\tran {(\cdot)}$,
$\Delta^+(\cdot) $, $\Delta^+(\cdot) $, $\difp {(\cdot)}$ and $\difn
{(\cdot)}$ also for vector-valued and tensor-valued functions,
intended as acting component-wise.

We will use the following properties of the difference
quotients, all proved in~\cite{hugo-petr-rose}. 
  Let $\bv\in W^{1,1}(\Omega)$ be such that $\spt \bv
  \subset\spt\xi$. Then 
\begin{equation}
\label{eq:commutation}
\begin{aligned}
  \nabla\difpm \bv &=\difpm{\nabla \bv }+\trap{(\partial_3 \bv
    )}\otimes\difpm{\nabla \Grenze}\,,
  \\
  \bD\difpm \bv &=\difpm{\bD \bv }+\trap{(\partial_3 \bv
    )}\otimess\difpm{\nabla \Grenze}\,,
  \\
  \diver\difpm \bv &=\difpm\diver \bv +\trapm{(\partial_3 \bv
    )}\difpm{\nabla \Grenze}\,,
  \\
  \nabla \bv _{\pm\tau} &= (\nabla \bv )_{\pm\tau} + \trapm{(\partial_3 \bv
    )}\difpm{\nabla \Grenze}\,,
\end{aligned}
\end{equation}
where $(\bv \otimes \bw)_{ij}:=v_iw_j$, $i,j=1,2,3$, and
$\bv\otimess \bw :=\frac 12 \big(\bv\otimes \bw + (\bv\otimes
\bw)^\top \big)$. 
%
Moreover,  we have also the following properties: 
  If $\spt g \subset\spt\xi$, then  there holds 
\begin{equation*}
  \trap{(\difn g )}=-\difp g \,,\quad \tran{(\difp g )}=-\difn g \,, \quad 
  \difn  g_\tau = - \difp g\,,
\end{equation*}
and
  if $\spt g\cup\spt f\subset\spt\xi$, then we have 
  \begin{equation*}
  \difpm (f g) = f_{\pm\tau} \,\difpm g + (\difpm f )\, g\,.
\end{equation*}
As for the classical difference quotients, $L^{q}$-uniform bounds (with
respect to $h>0$) for $\difp f$ imply that $\partial _\tau f$
belongs to $L^q(\spt\xi)$.
\begin{lemma}
  \label{lem:Dominic}
  If
  $ f \in W^{1,1}(\Omega)$, then we have for $\alpha=1,2$
  \begin{align}
    \label{eq:1}
    \difp f \to \td f=\partial _{\tau_\alpha}f :=\partial_\alpha f +\partial_\alpha
    \Grenze\, \partial_3 f  \qquad \text{ as } h\to 0\,,
  \end{align} 
  almost everywhere in $\spt\xi$, (cf.~\cite{mnr3}).
  If we
  define, for $0<h<R_P$
  \begin{equation*}
    \Omega_{P,h}=\left\{\bx\in \Omega_P\fdg x'\in B^2_{R_P-h}(0)\right\}\,,
  \end{equation*}
  and, if $f\in W^{1,q}_\loc(\setR^3)$, for $1\leq q<\infty$, then
  \begin{equation*}
    \int\limits_{\Omega_{P,h}}|d^+f|^q\,d\bx\leq c\int\limits_{\Omega_{P}}|\partial_\tau f|^q\,d\bx\,.
  \end{equation*}
 Moreover, if $d^{+}f\in
L^q(\Omega_{P,h_0})$, $1<q<\infty$, and if
  \begin{equation*}
    \exists\,c_1>0:\quad   \int\limits_{\Omega_{P,h_0}}|d^{+}f|^q\,d\bx\leq c_1\qquad
    \forall\,h_0\in(0,R_P)\text{ and } \forall\,h\in(0,h_0)\,,
  \end{equation*}
  then $\partial_{\tau}f\in L^q(\Omega_P)$ and
  \begin{equation*}
    \int\limits_{\Omega_{P}}|\partial_{\tau}f|^q\,d\bx\leq c_1\,.
  \end{equation*}
\end{lemma}
The following variants of formula of integration by parts will 
often be used.
\begin{lemma}
  \label{lem:TD3}
  Let $\spt g\cup\spt f\subset\spt\xi=\spt\xi_P$ and $0<h<\frac{R_P}{16}$. Then
  \begin{equation*}
    \intO f\tran g \, d\bx =\intO\trap f g\, d\bx\,.
  \end{equation*}
  Consequently, $\intO f\difp g \, d\bx= \intO(\difn f )g\, d\bx$\,.
  Moreover, if in addition $f$ and $g$ are smooth enough and at least
  one vanishes on $\partial\Omega$, then 
\begin{equation*} \intO f\td g \, d\bx= -\intO(\td f )g\,
    d\bx\,.
  \end{equation*}
\end{lemma}
\vspace*{-6mm}
\subsection{Regularity results with possible dependencies on $A$}
\vspace*{-2mm}
We start proving spatial regularity for the approximate problem. %
The estimates, which will be proved for first order derivatives of $\nabla\bue$ and
$\bF^{A}(\bD\bue)$ in this first step, are uniform with respect to
$A\geq1$:
\begin{description}[(ii)]
\item[(i)] in the interior of $\Omega$; 
\item[(ii)] for tangential derivatives near the boundary.
\end{description}
On the contrary, the estimates depend on $A$ in the
normal direction near  the boundary $\partial\Omega$. Nevertheless, this allows later on to use the equations point-wise and
to prove (in a different way) estimates independent of $A\geq1$ even near the
boundary, allowing then to pass to the limit with $A\to \infty$. 

By using the translation method we obtain the following results, which
will be proved below. 
\begin{proposition}  
  \label{prop:JMAA2017-1}
  Let the operator $\bS=\partial \pot$, derived from the potential
  $\potV$, have $(p,\delta)$-structure for some $p\in(1,2]$, and
  $\delta\in(0,\infty)$, with characteristics
  $(\gamma_{3},\gamma_{4},p)$. Let $\Omega\subset\setR^3$ be a bounded
  domain with $C^{2,1}$ boundary and assume that
  $\bff \in L^{p'}(\Omega)$.  Then, the unique weak solution
  $\bue\in W^{1,2}_{0}(\Omega)$ of the approximate
  problem~\eqref{eq:plasticity-app} satisfies
  \begin{align}
 \label{eq:est-app}
    \begin{aligned}
      \int\limits_{\Omega} \xi_0^2
      \abs{\nabla \bF^{A}(\bD\bue)}^2
      +
      \function^{A}\big( \xi_{0}^2|\nabla^{2}\bue|\big)+(\delta+A)^{p-2} \xi_{0}^2|\nabla^{2}\bue|^2\,d\bx
     &\le c_{0} \,,
      \\[3mm]
      \int\limits_{\Omega} \xi^2_P \abs{\td
        \bF^{A}(\bD\bue)}^2+
     \function^{A}\big(
     \xi_{P}^2|\td\nabla\bue|\big)+(\delta+A)^{p-2}
     \xi_{P}^2|\td\nabla\bue|^2\,d\bx\, 
      &
      \le c_{P}\,,
        \end{aligned}
  \end{align}
  where
  $c_{0}=c_{0} (\para, \|\bff\|_{p'},\norm{\xi_0}_{
    1,\infty},\gamma_3,\gamma_4,p)$, while the constant related to the
  neighborhood of $P$ is such that
  $c_{P}=c_{P}
  (\para,\|\bff\|_{p'},\norm{\xi_P}_{1,\infty},\norm{\Grenze_{P}}_{C^{2,1}},$
  $\gamma_3,\gamma_4,p)$.
  Here, $\xi_{0}(\bx)$ is a cut-off function with support in the interior of
  $\Omega$ and, for arbitrary $P\in \partial \Omega$, the tangential
  derivative 
  is defined locally in $\Omega_P$ by~\eqref{eq:1}. 
\end{proposition}
By using Proposition~\ref{prop:JMAA2017-1} and the ellipticity of $\bS^A$ we can write, for
a.e. $\bx\in \Omega$, the missing partial derivatives in the normal direction (which is
locally $\be_{3}$ after a rotation of coordinates) in terms of the tangential ones. By
employing the previous results we obtain estimates also for the partial derivatives in the
$\be_{3}$-direction, but with a critical dependence on the approximation parameter $A$.
\begin{proposition}  
  \label{prop:JMAA2017-2}
  Under the assumptions of Proposition \ref{prop:JMAA2017-1} there
  exists a constant $C_1>0$ such
  that, 
  provided in the local description of the boundary there holds
  $r_P<C_1$ in $(b3)$, where $\xi_{P}(\bx)$ is a cut-off function with
  support in $\Omega_P$, there holds
  \begin{equation*}
    \begin{aligned}
\int\limits_{\Omega} 
      \xi^2_{P} \abs{\partial _3 \bF^{A}(\bD\bue)}^2
      +\function^{A}\big(\xi^2_{P}    \abs{\partial_3 \bD\bue}\big)\,d\bx
      \le C_{A}\,,
    \end{aligned}
  \end{equation*}
  where
  $C_{A}=C_{A}(\delta,\|\bff\|_{p'},\norm{\xi_P}_{1,\infty},\norm{\Grenze_{P}}_{C^{2,1}},\gamma{}_3,\gamma{}_4,p, 
  A)$.
\end{proposition}
Before starting the proof of these two propositions, we
generalize~\cite[Lemma~3.11]{br-reg-shearthin}, originally proved for
$\phi=\function$ (with $1<p\leq2$) to
$\phi=\function^{A}$. The main properties used are: convexity of
$\function^{A}$, that $(\function^{A})''$ is non-increasing, and the
equivalence properties from Lemma~\ref{lem:eq12}.
\begin{lemma}
  \label{lem:p-est}
Let $p \in (1,2]$ and
  $\delta\ge 0$. Then, for $\xi $ and $\Grenze$ as above and
 for any   $\bv \in W^{1,2}_{0}(\Omega)$ we have
  \begin{equation*}
    \begin{aligned}
     \intO \function^{A}\big(\xi \abs{\nabla \difp\bv}\big) + \function^{A}\big(\xi
      \abs{ \difp \nabla \bv}\big)\, d\bx &\le c\,\intO \xi^2 \bigabs{
        \difp \bF^{A} (\bD\bv) }^2 \, d\bx
      \\
      &\quad + c( \norm{\xi}_{1,\infty},
      \norm{\Grenze}_{C^{1,1}}) \hspace*{-3mm}\int\limits_{\Omega\cap \spt                   
        \xi}\hspace*{-3mm} \function^{A} \big (\abs{ \nabla\bv }\big ) \,
      d\bx\\[-4mm]
    \end{aligned}
  \end{equation*}
  with constants not depending on $\delta$ and $A$.
\end{lemma}
\begin{proof}
  The proof is carried out by adapting that
  of~\cite[Lem.~3.11]{br-reg-shearthin}. 
%
%
%
First, we use the following identity
  \begin{equation*}
    \xi\, \nabla \difp \bv = \nabla (\xi\, \difp \bv )
    -\nabla\xi\otimes \difp \bv\,,
  \end{equation*}
 and consequently we get, by using~\eqref{eq:lambda}, that 
 \begin{equation*}
   \begin{aligned}
     \intO \function^{A} (\xi \abs{\nabla \difp \bv}) \, d\bx &\le c \intO \function^{A}
     (\abs{\bD (\xi \,\difp \bv)}) \, d\bx + c(\norm{g}_{C^{0,1}},\norm{\xi}_{1,\infty}) \hspace*{-3mm}
     \int\limits_{\Omega\cap \spt \xi}\hspace*{-3mm}\function^{A} (\abs{\nabla \bv}) \, d\bx \,,
   \end{aligned}
 \end{equation*}
 where we also used Korn's inequality for N-functions
 (cf.~\cite[Thm.~6.10]{john}), with a constant independent of
 $A\geq1$,  and the following inequality (cf.~\cite[Sec.~3.2]{br-multiple-approx})
 \begin{equation}
   \label{eq:diff}
   \int\limits_{\Omega\cap \spt \xi} \hspace*{-3mm}\function^{A}(\abs{\difpm\bv}) \, d\bx \le 
   c\, \hspace*{-3mm}\int\limits_{\Omega\cap \spt \xi}
   \hspace*{-4mm}\function^{A}(\abs{\nabla\bv}) \, d\bx \,. \\[-3mm]
 \end{equation}
 Using the identities
\begin{align*}
  \bD(\xi\, \difp\bv) &= \xi\, \bD(\difp \bv) + \nabla \xi \otimess
  \difp \bv
     = \xi\, \difp \bD\bv + \trap{(\partial _3 \bv)} \otimess
  \difp\nabla\Grenze + \nabla \xi \otimess \difp \bv\,,
\end{align*}
the properties of $\function^{A}, \xi, \Grenze$, and~\eqref{eq:diff}, we obtain
\begin{equation}
  \label{eq:p-est2}
  \begin{aligned}
    \intO \hspace*{-1mm} \function^{A} (\xi \abs{\nabla \difp \bv}) \, d\bx
    &\le c\hspace*{-1mm}\intO \hspace*{-1mm}\function^{A} (\xi \abs{\difp
      \bD\bv}) \, d\bx
     \\
     &\quad
    + c(\norm{\xi}_{1,\infty}, \norm{\Grenze}_{C^{1,1}})
    \hspace*{-3mm} \int\limits_{\Omega\cap \spt \xi}\hspace*{-3mm}\function^{A}
    (\abs{\nabla \bv}) \, d\bx \,.
  \end{aligned}
\end{equation}
We focus on the first term on the right-hand side
of~\eqref{eq:p-est2}. Using a change of
shift as in Lemma~\ref{lem:change2} yields
\begin{align}\label{eq:shi}
    \function^{A} (\xi \abs{\difp \bD\bv}) &\le c\big (\function^{A}_{\abs{\bD\bv}
    +\abs{\Delta^{+} \bD\bv}} (\xi \abs{\difp \bD\bv}) +\function^{A}
                                             (\abs{\bD\bv} + 
  \abs{\Delta^+ \bD\bv}) \big )\,.
\end{align}
\enlargethispage{5mm}
By the fact that $\function^{A}$ is balanced with characteristics depending
only on $p$ (cf.~Lemma~\ref{lem:function}, Lemma \ref{lem:eq12}) we get
\begin{align*}
  \function^{A}_{\abs{\bD\bv}
  +\abs{\Delta^+ \bD\bv}} (\xi \abs{\difp \bD\bv})
   &
     \sim   (\function^{A}_{\abs{\bD\bv} +\abs{\Delta^+ \bD\bv}})' (\xi
     \abs{\difp \bD\bv})\
     \xi \abs{\difp \bD\bv}
  \\
   &=\frac{(\function^{A})' (\abs{\bD\bv} +\abs{\Delta^+ \bD\bv}+\xi
     \abs{\difp \bD\bv})}{\abs{\bD\bv} +\abs{\Delta^+ \bD\bv}+\xi
     \abs{\difp \bD\bv}}
     \ \xi^{2} \abs{\difp \bD\bv}^{2}
  \\
   &=a^{A} (\abs{\bD\bv} +\abs{\Delta^+ \bD\bv}+\xi
     \abs{\difp \bD\bv})     \ \xi^{2} \abs{\difp \bD\bv}^{2}.
\end{align*}
Next, since $a^{A}$ is non-increasing for $p \in (1,2]$ (see
Lemma~\ref{lem:UA}) we get 
\begin{align*}
  \function^{A}_{\abs{\bD\bv}
  +\abs{\Delta^+ \bD\bv}} (\xi \abs{\difp \bD\bv})
  &\leq a^{A} (\abs{\bD\bv} +\abs{\Delta^+ \bD\bv})     \ \xi^{2}
    \abs{\difp \bD\bv}^{2}
 %
  \\
  &\sim \xi^{2} \abs{\difp \bF^{A}(\bD\bv)}^{2}\,.\\[-4mm]
\end{align*}
 Inserting this into~\eqref{eq:shi} we obtain from \eqref{eq:p-est2}
 that 
 \begin{equation}
    \label{eq:p-est3}
   \begin{aligned}
     &\intO \function^{A} (\xi \abs{ \nabla \difp \bv}) \, d\bx
     \\
     &\le c \intO
    \xi^{2} \bigabs{\difp \bF^{A}(\bD\bv )}^2 \, d\bx
     + c(\norm{\xi}_{1,\infty}, \norm{\Grenze}_{C^{1,1}}) \hspace*{-3mm}
     \int\limits_{\Omega\cap \spt 
       \xi}\hspace*{-3mm}\function^{A} (\abs{\nabla \bv}) \, d\bx \,. \\[-4mm]
   \end{aligned}
 \end{equation}
Next, we observe that \eqref{eq:commutation} and \eqref{eq:lambda} yield 
\begin{equation*}
  \begin{aligned}
    \intO \function^{A} (\xi \abs{ \difp \nabla \bv}) \,d\bx&\leq \intO \function^{A} (\xi
    \abs{ \nabla \difp\bv}) \,d\bx+ \intO \function^{A} (\xi
    \abs{\partial_3\bv}\abs{\difp\nabla\Grenze}) \,d\bx
    \\
    &\leq \intO \function^{A} (\xi \abs{ \nabla \difp\bv}) \,d\bx+ c(\norm{\xi}_{\infty},
    \norm{\Grenze}_{C^{1,1}}) \hspace*{-3mm}\int\limits_{\Omega\cap \spt \xi}\hspace*{-3mm}\function^{A} (\abs{\nabla
      \bv}) \, d\bx \,. \\[-4mm]
  \end{aligned}
\end{equation*}
This shows that also the term $ \int_{\Omega} \function^{A} (\xi \abs{ \difp \nabla \bv}) \,d\bx$ can
be estimated by the right-hand side of~\eqref{eq:p-est3}, ending the proof.
\end{proof}
%
We can now proceed with the proof of regularity in the tangential
directions and  in the interior.
%
\begin{proof}[of Proposition~\ref{prop:JMAA2017-1}]
  We obtain estimates for tangential derivatives by considering
  limits of increments in the tangential directions
  cf.~\cite{br-reg-shearthin,br-multiple-approx}. Fix
  $P\in \partial \Omega$ and use in $\O_P$
  \begin{equation*}
    \bw=\difn{(\xi^2\difp(\bue |_{\frac 12 {\Omega}_P}))}\,,
  \end{equation*}
  where $\xi:=\xi_P$, $\Grenze:=\Grenze_P$, and
  $h\in(0,\frac{R_{P}}{16})$, as a test function in the weak
  formulation~\eqref{eq:weak-app} of
  Problem~\eqref{eq:plasticity-app}. This yields
\enlargethispage{5mm}
  \begin{equation}
  \label{eq:dtTx2}
  \begin{aligned}
    &\intO \xi^2\difp{\bS^{A}(\bD\bue)}\cdot \difp \bD\bue \, d\bx
    \\
    &=-\intO \bS^{A}(\bD\bue)\cdot\big(\xi ^2 \difp \partial _3 \bue
    -(\xi_{-\tau } \difn\xi +\xi \difn\xi) \partial_3\bue \big)
    \otimess\difn\nabla\Grenze\, d\bx
    \\
    &\quad-\intO \bS^{A}(\bD\bue)\cdot \xi^2 \trap{(\partial _3\bue)}\otimess
    \difn{\difp{\nabla\Grenze}} - \bS^{A}(\bD\bue)\cdot\difn{\big(2\xi \nabla \xi
      \otimess \difp\bue\big)}\, d\bx
    \\
    &\quad+\intO \bS^{A}(\trap{(\bD\bue)})\cdot \big(2 \xi \partial_3\xi \difp\bue +
    \xi^2 \difp\partial_3\bue \big)\otimess\difp\nabla\Grenze \, d\bx
    \\
    &\quad+\intO\ff\cdot\difn(\xi^2 \difp \bue)\, d\bx=:\sum_{j=1}^{8} I_j\,.
  \end{aligned}
\end{equation}
The properties of $\bS^A$, Proposition~\ref{prop:SA-ham}, and
Lemma~\ref{lem:p-est} imply the following estimate
\begin{equation*}
  \begin{aligned}
    \label{eq:odhadT2}
    &\intO \xi^2 
    \bigabs{ \difp \bF^{A} (\bD\bue) }^2
    + \function^{A}\big(\xi
      \abs{ \difp \nabla \bu^{A}}\big)\, d\bx     
    \\
    &\leq   c \intO\xi ^2 \difp{\bS^{A}(\bD\bue)}\cdot
    \difp \bD\bue\,d\bx + c( \norm{\xi}_{1,\infty}, 
      \norm{\Grenze}_{C^{1,1}}) \hspace*{-3mm}  \int\limits_{\Omega\cap \spt
        \xi}\hspace*{-3mm} \function^{A} \big (\abs{ \nabla\bue }\big ) \,
      d\bx\,.
  \end{aligned}
\end{equation*}
The terms $I_1$--$I_7$ in~\eqref{eq:dtTx2} are estimated exactly as
in~\cite[(3.17)--(3.22)]{br-reg-shearthin}, while $I_8$ is estimated as
the term $I_{15}$ in~\cite[(4.20)]{br-reg-shearthin}. Thus, we get, by
using also Corollary~\ref{cor:UA}, 
\begin{align*}
  & \intO (\delta+A)^{p-2}\xi^2\bigabs{\difp\nabla
    \bue }^2 \! +\!\xi^2 \bigabs{ \difp \bF^{A} (\bD\bue) }^2
    \!+\! \function^{A} (\xi 
    \abs{\difp\nabla \bue}) 
    \notag
  \\
  & \le c(\|\ff\|_{p'}, \norm{\xi}_{2,\infty},\norm{\Grenze}_{C^{2,1}},
    \delta )\,. 
\end{align*}
This proves the second estimate in~\eqref{eq:est-app} by
Lemma~\ref{lem:Dominic}, since the constant on the right-hand side does not depend on
$h>0$.

The first estimate in~\eqref{eq:est-app} is proved in the same way
with many simplifications, since in the interior one can consider
directly standard translations in all the coordinate directions. 
\end{proof}
For the proof of Proposition~\ref{prop:JMAA2017-2}, the following
observation will be crucial.
\begin{remark}
  \label{rem:point-wise}
  The obtained estimate $\eqref{eq:est-app}_1$, Proposition
  \ref{prop:pFA}, and Lemma \ref{lem:UA} imply that $\bue \!\in\!
  W^{2,2}_{\loc }(\Omega)$ (with estimates depending on $A$)
  and that the system~\eqref{eq:plasticity-app} is {well-defined} point-wise
  a.e.~in $\Omega$.
\end{remark}
\begin{proof}[of Proposition~\ref{prop:JMAA2017-2}]
  To estimate the derivatives in the $\be_{3}$-direction we use
  equation~\eqref{eq:plasticity-app} point-wise a.e.~in $\Omega$,
  which is justified by Remark~\ref{rem:point-wise}.  Denoting, for
  $\alpha,\gamma=1,2$,
  $A_{\alpha \gamma}:=\partial_{\gamma 3}S^{A}_{\alpha 3}(\bD\bue)$,
  $\mathfrak {b}_\gamma:=\partial_3 D_{\gamma 3}\bue$,
  and\footnote{Recall that we use the summation convention over
    repeated Greek lower-case letters from $1$ to $2$. }
  $\mathfrak {f}_\alpha :=f_\alpha+\partial_{3 3}S^{A}_{\alpha
    3}(\bD\bue)\partial_3 D_{3 3}\bue +\partial_{\gamma
    \sigma}S^{A}_{\alpha 3}(\bD\bue)\partial_3 D_{\gamma \sigma}\bue+
  \sum_{k,l=1}^3\partial_{k l}S^{A}_{\alpha
    \beta}(\bD\bue)\partial_\beta D_{k l}\bue$, we can re-write the
  first two equations in~\eqref{eq:plasticity-app} as follows:
  \begin{equation*}
    \label{eq:linear_system}
    -2  A_{\alpha\gamma}\mathfrak b_\gamma=\mathfrak{f}_\alpha\qquad\textrm
    {a.e.\ in }\Omega\,.
  \end{equation*}
  We employ this equality separately on each $\Omega_{P}$ in order to use the notion
  of tangential derivative. By straightforward manipulations
  (cf.~\cite[Sections 3.2 and 4.2]{br-reg-shearthin}) we get 
  a.e.~in $\Omega_{P}$
  \begin{equation*}
    \begin{aligned}
      & \MC^{A}(|\bD\bue|)\,| {\boldsymbol {\mathfrak b}}|
      \\
      &\leq c
      \left(|\bff|+|\bff|\|\nabla\Grenze\|_{\infty}+\MC^{A}(|\bD\bue|)
        \,\big(|\partial_\tau\nabla
          \bue|+\|\nabla\Grenze\|_{\infty}|\nabla^2 \bue|\big)\right)\,.
    \end{aligned}
  \end{equation*}
  Note that we can deduce from this inequality information about
  $\tilde {\mathfrak b}_\gamma:=\partial ^2_{33}\bue_{\gamma}$,
  because
  $ |\boldsymbol { \mathfrak b}|\geq2|\tilde{\boldsymbol{\mathfrak
      b}}|-|\partial_{\tau}\nabla\bue|-\|\nabla\Grenze\|_{\infty}|\nabla
  ^{2}\bue|$. Adding on both sides, for $\alpha=1,2$ and $i,k=1,2,3$
  the term
  \begin{equation*}
\MC^{A}(|\bD\bue|)\,\big(|\partial_\alpha\partial_i
    u^{A}_{k}|+|\partial_{33}^2 u^{A}_{3}|\big)\,,
  \end{equation*} 
  we finally arrive,  a.e.~in $\Omega_{P}$ at the inequality
  \begin{equation*}
    \begin{aligned}
      & \MC^{A}(|\bD\bue|)| {\nabla ^2\bue }|
      \\
      &\leq c \left(|\bff|+|\bff|\|\nabla
        g\|_{\infty}+\MC^{A}(|\bD\bue|)\,\big(|\partial_\tau\nabla
          \bue|+\|\nabla\Grenze\|_{\infty}|\nabla^2 \bue|\big)\right)\,,
    \end{aligned}
  \end{equation*}
  where, due to the results proved in Section~\ref{sec:section-2}, the
  constant $c$ only depends on the characteristics of $\bS$.  Next, we
  can choose the open sets $\Omega_{P}$ in such a way that
  $\|\nabla\Grenze_{P}(x')\|_{{\infty},\Omega_{P}}$ is small enough, so
  that we can absorb the last term from the right-hand side, which
  yields \\[-2mm]
  \begin{equation*}
    \label{eq:estimate_normal_final}
    \begin{aligned}
      &  \MC^{A}(|\bD\bue|)\,|\nabla ^{2}\bue|
       \leq c\left(|\bff|+
         \MC^{A}(|\bD\bue|)\,|\partial_{\tau}\nabla\bue|\right)
       \quad\textrm{ a.e. in }\Omega_P\,,
    \end{aligned}
  \end{equation*}
  where again the constant $c$ only depends on the characteristics of
  $\bS$.  Dividing both sides by the quantity
  $\sqrt{ \MC^{A}(|\bD\bue|)}\not=0$ (which is non zero by the fact
  that $\delta>0$ and the properties of $\function^{A}$) and raising
  the result to the power $2$ we get a.e.~in $\Omega_{P}$
  \begin{equation}
    \label{eq:6}
    \begin{aligned}
     \MC^{A}(|\bD\bue|)|\nabla^2\bue|^2&\leq
     c\, \frac{|\bff|^2}{ \MC^{A}(|\bD\bue|)}
      +c\, 
      \MC^{A}(|\bD\bue|)|\partial_\tau\nabla\bue|^2\,.
    \end{aligned}
  \end{equation}
  Note that both sides are finite a.e. and, for the moment, we know that
  the left-hand side belongs at least to $L^{1}_{loc}(\Omega_{P})$.

  Concerning the first term on the right-hand side, we note that
  Lemma~\ref{lem:UA} and the definition of $\MC^A$ imply \\[-2mm]
  \begin{equation*}
    \frac{1}{\MC^{A}(t)}\leq     \frac{1}{(p-1)
    }\frac{1}{\MC(t)}=\frac{1}{(p-1) }\frac{1}{(\delta +t)^{p-2}}\,.
  \end{equation*}
  Using this estimate  and 
  H\"older inequality 
  we get, with a constant $c$ independent on~$A$,
  \begin{equation}
  \label{eq:estimate-f}
  \begin{aligned}
    \intO\frac{|\bff|^2}{ \MC^{A}(|\bD\bue|)}\,d\bx&\leq
    c(p)\|\bff\|_{p'}^{2}\|\delta+|\bD\bue|\|_{p}^{2-p}
    \\[-3mm]
    &\leq c\big(\|\bff\|_{p'}^{p'}+\delta^{p}+\|\bF(\bD\bue)\|_{2}^{2}\big)\,.
  \end{aligned}
\end{equation}
%
  For the second term on the right-hand side of \eqref{eq:6} we use that,
  in view of Lemma~\ref{lem:UA}, there holds:
  \begin{equation*}
    \begin{aligned}
      \intO \xi_P^2
      \MC^{A}(|\bD\bue|)|\partial_\tau\nabla\bue|^2\,d\bx&\leq
      \delta^{p-2}\intO \xi_P^2|\partial_\tau\nabla\bue|^2\,d\bx
      \\
      &=   \frac{\delta^{p-2}}{(\delta+A)^{p-2}}\ (\delta+A)^{p-2}\intO \xi_P^2|\partial_\tau\nabla\bue|^2\,d\bx
      \\[-1mm]
      & \leq    \left(1+\frac{A}{\delta}\right)^{2-p}\, c_{P} \,,
    \end{aligned}
  \end{equation*}
  where the final estimate follows from the the already proved results on
  tangential derivatives in Proposition \ref{prop:JMAA2017-1}.

Hence, multiplying \eqref{eq:6} by $\xi_{P}^2$ and integrating over 
 the proper sub-domain
  \begin{equation*}
\Omega_{P,\epsilon}:=\left\{\bx\in \Omega_{P}\fdg
  \Grenze_{P}+\epsilon<x_{3}<\Grenze_{P}+R_{P}',\text{ for }
  0<\epsilon<R_{P}'\right\}\,, 
  \end{equation*}
we get, also using \eqref{eq:main-apriori-estimate2} and \eqref{eq:f},
  \begin{equation*}
    \begin{aligned}
      &\int\limits_{\Omega_{P,\epsilon}} \xi_P^2
      \MC^{A}(|\bD\bue|)|\nabla^2\bue|^2\,d\bx
      \\
      &\leq c\int\limits_{\Omega_{P,\epsilon}}\frac{|\bff|^2}{
        \MC^{A}(|\bD\bue|)}\,d\bx +c \int\limits_{\Omega_{P,\epsilon}}
      \xi_P^2 \MC^{A}(|\bD\bue|)|\partial_\tau\nabla\bue|^2\,d\bx
      \\
      &\leq c\intO\frac{|\bff|^2}{ \MC^{A}(|\bD\bue|)}\,d\bx +c \intO
      \xi_P^2 \MC^{A}(|\bD\bue|)|\partial_\tau\nabla\bue|^2\,d\bx
      \\
      &\leq C\, \big (\delta^{p} +\norm{\bff}_{p'}^{p'} + (1+A
      \delta^{-1})^{2-p} \, \big)\,.
    \end{aligned}
  \end{equation*}
  Since this estimate is independent of $\epsilon>0$, the above
  inequality shows, by monotone convergence, that also
  $\int_{\Omega} \xi_P^2 \MC^{A}(|\bD\bue|)|\nabla^2\bue|^2\,d\bx\leq
  C(\norm{\bff}_{p'},\delta,\delta^{-1},A)$, ending the proof.
\end{proof}
\begin{remark}
  The reader should notice that the dependence on $A$ is mainly due to the
  fact that we have a stress tensor depending on the symmetric
  gradient. To use Korn inequality in Lemma~\ref{lem:p-est} we have to
  pay the price of estimates depending on $A$. In the case of a stress
  tensor depending on the full gradient, this step can be skipped (see
  the results in \cite{BG2016} where the considered problem has an
  additional term with $2$-structure and the $A$-approximation is not
  needed).
\end{remark}
By collecting the results of the
Propositions~\ref{prop:JMAA2017-1} and \ref{prop:JMAA2017-2} we get the
following result.
\begin{proposition}
  \label{prop:regularity-A}
  Let the operator $\bS=\partial \pot$, derived from the potential
  $\pot$, have $(p,\delta)$-structure for some $p \in (1,2]$, and
  $\delta \in (0,\infty)$.
Let $\O\subset\setR^3$ be a bounded domain with $C^{2,1}$
  boundary and let $\ff\in L^{p'}(\O)$. Then, the unique weak solution
  $\bue\in W^{1,2}_{0}(\Omega)$ of  problem~\eqref{eq:plasticity-app}
  satisfies \\[-2mm]
  \begin{equation*}
    \begin{aligned}
      \intO   \abs{\nabla \bF^{A}(\bD\bue)}^2\,d\bx
      &\le c(A,\delta^{-1}) \,, \\[-2mm]
    \end{aligned}
  \end{equation*}  
  where $c$ depends also on the characteristics of $\bS$, $\delta$, $\|\ff\|_{p'}$, $\abs{\Omega}$, and the
  $C^{2,1}$-norms of the local description of $\partial \Omega$.  In
  particular, the above estimate implies that
  ${\bue\in W^{2,2}(\Omega).  }$
\end{proposition}
\begin{proof}
  The proof is simply obtained by observing that $\overline{\Omega}$
  is a compact set. After having fixed all $\Omega_{P}$ small enough
  (depending on $P$) to perform the calculations leading to
  Proposition~\ref{prop:JMAA2017-2}, we can extract a finite covering
  of sets $\Omega_{P}$ and consequently prove the
  uniform nature of estimates in terms of $P\in \partial\Omega$. To
  show that $\bue\in W^{2,2}(\Omega) $ we use Proposition
  \ref{prop:pFA} and Lemma   \ref{lem:UA}. 
\end{proof}
%
%
\section{Estimates uniform with respect to $A$ for the
  solutions of the approximate problem}
\label{sec:4}%
We now sketch the proof of the estimate of $\nabla\bF^{A}(\bD\bue)$,
which is independent of~$A$. Moreover, we also improve the
$\delta$-dependence of the estimates.  We adapt to the
$A$-approximation the same procedure already used
in~\cite{br-plasticity}, obtaining the steady counterpart to the case
$1<p\leq2$ of the results proved in~\cite[Sec.~3]{br-multiple-approx}.
\begin{proposition}
  \label{prop:main}
  Let the same hypotheses as in Theorem~\ref{thm:MT} be satisfied with
  $\delta >0$ and let the local description $\Grenze_{P}$ of the
  boundary and the localization function $\xi_P$ satisfy
  $(b1)$--\,$(b3)$ and $(\ell 1)$ (cf.~Section~\ref{sec:bdr}). Then,
  there exists a constant $C_2>0$ such that the regular solution
  $\bue\in W^{1,2}_{0}(\Omega)\cap W^{2,2}(\Omega)$ of the approximate
  problem~\eqref{eq:plasticity-app} satisfies for every
  $P\in \partial \Omega$
  \begin{equation*}
    \int\limits_\Omega 
      \xi^2_P |\partial_3\bF^{A}(\bD\bue)|^2\,d\bx
     \leq C\,, 
  \end{equation*}
  provided $r_P<C_2$ in $(b3)$, with $C$ depending on the
  characteristics of $\bS$, $\delta$, 
  $\norm{\bff}_{p'}$, $\norm{\xi_P}_{
    1,\infty},\norm{\Grenze_{P}}_{C^{2,1}}$, and $C_2$.
\end{proposition}
\begin{proof}
  We will not give the full proof of this result, since it is very similar to that of
 \cite[Prop.~3.2]{br-plasticity}. For the reader's convenience we just explain the
  main steps.

Fix an arbitrary point
  $P\in \partial \Omega$ and a local description $\Grenze=\Grenze_{P}$
  of the boundary and the localization function $\xi=\xi_P$ as
  before. 
  Proposition~\ref{prop:pFA} yields that there exists a
  constant $C_0$, depending only on the characteristics of
  $\bS$ 
  such that
  \begin{equation*}
   \frac{1}{C_0}| \partial_3\bF^{A}(\bD\bue)|^2\leq
   \mathbb{P}_3^{A}(\bD\bue)\qquad \text{a.e.  in }\Omega\,.
 \end{equation*}
 We now work directly with $\mathbb{P}_3^{A}(\bD\bue)$ 
 to deduce estimates for $|\partial_3\bF^{A}(\bD\bue)|^{2}$. Note
 that,  since $\bue$ is a regular solution of \eqref{eq:plasticity-app}, all
 calculations are justified.
%
%
 Thus, using the definition of  $\mathbb{P}_3^{A}(\bD\bue)$ and the
 symmetries of $\bS^A$ and $\bD\bue$, we obtain 
  \begin{align}
    \label{eq:j123}
    &\frac{1}{C_0}
      \int\limits_\Omega 
      \xi^2 |\partial_3\bF^{A}(\bD\bue)|^2\,d\bx 
    \\
    & \leq \int\limits_\Omega\xi^2 
      \partial_3 \bS^{A}_{\alpha\beta} (\bD\bue) \,\partial_3
      D_{\alpha\beta}\bue\,d\bx  +
      \int\limits_\Omega\xi^2
      \partial_3\bS^{A}_{3\alpha }(\bD\bue)\,\partial _\alpha
      D_{33}\bue\,d\bx \notag
    \\
    &\quad +\int\limits_\Omega \sum_{j=1}^3
      \xi^2\partial_3  \bS^{A}_{j3}(\bD\bue)\, \partial_3^2
      u^{N}_{j}\,d\bx \notag 
    \\
    & =:\mathcal{J}_{1}+\mathcal{J}_{2}+\mathcal{J}_{3}\,.\notag
  \end{align}
  The most critical term is $\mathcal{J}_{1}$ which is estimated, for any
  $\lambda>0$, as follows 
\begin{align*} 
    |\mathcal{J}_{1}| \leq &\, \param
    \int\limits_\Omega\xi^2
    |\partial_3\bF^{A}(\bD\bue)|^2 \,d\bx 
+c_{\param^{-1}}\, \big (1+\|\nabla \Grenze\|_\infty^2\big )
    \sum_{\beta=1}^2 \int\limits_\Omega\xi^2
    |\partial_\beta\bF^{A}(\bD\bue)|^2 \,d\bx 
    \\
    &\quad+ \int\limits_\Omega\xi^2 
      |\partial_3 \bS^{A} (\bD\bue) |\, |\nabla
      ^2\Grenze| \, |\bD\bue|\,d\bx  
 + \bigg |\int\limits_\Omega\xi^2 
      \partial_3 \bS^{A}_{\alpha\beta} (\bD\bue) \,\partial_\alpha
     \partial _{\tau_\beta} u_3^{A}\,d\bx  \,\bigg|\,.
\end{align*}
In the last but one term 
we multiply and divide by
$\sqrt{\MC^{A}(|\bD\bue|)}$, use Proposition \ref{prop:pFA}, Young
inequality, and
$\MC^{A} (|\bD\bue|) |\bD\bue|^2 \sim |\bF^{A}(|\bD\bue|)|^2$
(cf.~Pro\-position \ref{prop:SA-ham}), yielding that it is estimated by
\begin{equation*}
  \param
    \int\limits_\Omega\xi^2
    |\partial_3\bF^{A}(\bD\bue)|^2 \,d\bx +c_{\param^{-1}} \, \|\nabla^2 \Grenze\|^2_\infty
  \int\limits_\Omega|\bF^{A}(\bD\bue)|^2 
\, d\bx \,.
\end{equation*}
To handle the last term in the above estimate of $\mathcal J_1$
we perform a crucial partial
integration. This avoids to have terms with the quantity
$\partial_3 \bS^{A} (\bD\bue)$ which cannot be estimated in terms of
tangential derivatives. Let us explain the main idea beyond this
step. Observe that, by neglecting the localization $\xi$, integration by
parts gives
\begin{align*}
  \int\limits_\Omega
      \partial_3 \bS^{A}_{\alpha\beta} (\bD\bue) \,\partial_\alpha
     \partial _{\tau_\beta} u_3^A\,d\bx  &=\int\limits_\Omega 
      \partial_\alpha \bS^{A}_{\alpha\beta} (\bD\bue) \,\partial_3
                                               \partial _{\tau_\beta} u_3^{A}\,d\bx 
  \\
  &=\int\limits_\Omega 
      \partial_\alpha \bS^{A}_{\alpha\beta} (\bD\bue) \,
      \partial _{\tau_\beta} D_{33} \bue\,d\bx \,.
\end{align*}
We next multiply and divide the integrand on the right-hand side by $\sqrt{\MC^{A}(|\bD\bue|)}$, use
Proposition~\ref{prop:pFA}, Young inequality, and the definition of the tangential
derivatives, yielding that
\begin{align*}
  \begin{aligned}
&\Big|\,\int\limits_\Omega 
      \partial_\alpha \bS^{A}_{\alpha\beta} (\bD\bue) \,
      \partial _{\tau_\beta} D_{33} \bue\,d\bx\,\Big|
     \\[-2mm]
      &\quad\leq c \sum_{\alpha=1}^2\int\limits_\Omega
    |\partial_\alpha\bF^{A}(\bD\bue)|^2 \,d\bx +c \sum_{\beta=1}^2
    \int\limits_\Omega|\partial _{\tau
      _\beta}\bF^{A}(\bD\bue)|^2 \, d\bx 
    \\[-1mm]
    &\quad \le c \sum_{\alpha=1}^2
    \int\limits_\Omega|\partial _{\tau
      _\alpha}\bF^{A}(\bD\bue)|^2 \, d\bx +c \,\|\nabla
    \Grenze\|_\infty ^2\int\limits_\Omega 
    |\partial_3\bF^{A}(\bD\bue)|^2 \,d\bx \,. \hspace*{-5mm}
  \end{aligned}
\end{align*}
The presence of the localization leads to several additional lower order terms, which all
can be easily handled as in \cite{br-plasticity}.
To
  treat $\mathcal J_2$ we multiply and divide by
  $\sqrt{\MC^{N}(|\bD\bue|)}$, use Proposition \ref{prop:pFA} and
  Young inequality, to show that, for any given $\lambda>0$, it holds
\begin{equation*}
  \begin{aligned}
    |\mathcal{J}_{2}| \leq &\,\param
    \int\limits_{0}^{t}\int\limits_\Omega\xi^2
    |\partial_3\bF^A(\bD\bue)|^2 \,d\bx  +c_{\param^{-1}}\,
    \sum_{\beta=1}^2 \int\limits_{0}^{t}\int\limits_\Omega\xi^2
    |\partial_\beta\bF^A(\bD\bue)|^2 \,d\bx \,,
  \end{aligned}
\end{equation*}
for some constant $c_{\param^{-1}} $ depending  on
$\param^{-1}$. To handle the term $\mathcal J_3$ we use the equation
\eqref{eq:plasticity-app}. All terms are handled exactly as in
\cite[Prop.~3.2]{br-plasticity}, and thus, we skip the details here. 
All together we arrive at the following
estimate
\begin{equation*}
  \begin{aligned}
    |\mathcal{J}_{1}|+  |\mathcal{J}_{2}| +|\mathcal{J}_{3}|
    &\leq \,\big (\param  
    +c_{\param^{-1}}\,\|\nabla \Grenze\|_\infty^2\big ) \int\limits_\Omega\xi^2
    |\partial_3\bF^{A}(\bD\bue)|^2 \,d\bx 
    \\
    &\ +c_{\param^{-1}}\sum_{\beta=1}^2 \int\limits_\Omega\xi^2
    |\partial_{\tau_\beta}\bF^{A}(\bD\bue)|^2 \,d\bx 
    \\
    &\  +c_{\param^{-1}} \big(1+\|\nabla\xi\|_\infty^2\big
    ) \!\int\limits_\Omega\!|\bF^{A}(|\bD\bue|)|^2\, d\bx
+    c_{\param^{-1}}\!
       \int\limits_{\Omega}\!\frac{|\bff|^{2}}{\MC^{A}(|\bD\bue|)}\,d\bx\,.
  \end{aligned}
\end{equation*}
Now we first choose $\lambda>0$ smaller than $(4 C_{0})^{-1}$ and then we choose the
covering of the boundary $\partial \Omega$ such that $c_{\param^{-1}}\,\|\nabla
\Grenze\|_\infty^2\leq(4 C_{0})^{-1}$, in order to absorb in the left-hand side of
\eqref{eq:j123} the term involving $\partial_3\bF^{A}(\bD\bue)$.  By using the estimate
\eqref{eq:estimate-f} already proved for the term with the external force, we get
\begin{equation*}
  \begin{aligned}
   \int\limits_\Omega 
\xi^2|\partial_3\bF^{A}(\bD\bue)|^2\,d\bx 
    &\leq    c 
    \sum_{\beta=1}^2     \int\limits_\Omega\xi^2|\partial_{\tau_{\beta}}\bF^{A}(\bD\bue)|^{2}\,d\bx
    \\
    &\quad
    +c\intO|\bF^{A}(\bD\bue)|^2\,d\bx
    +c\big(\delta^{p}+\norm{\bff}_{p'}^{p'}\big)\,
    \end{aligned}
\end{equation*}
with constants depending only on the characteristics of $\bS$, $\norm{g}_{C^{2,1}}$, and
$\norm{\xi}_{1,\infty}$. The uniform estimates~\eqref{eq:main-apriori-estimate2},
\eqref{eq:est-app} for the right-hand side allows us to end the proof of
Proposition~\ref{prop:main}.
\end{proof}
\begin{proposition}
  \label{prop:regularity-A-ind}
  Let the operator $\bS=\partial \pot$, derived from the potential
  $\pot$, have $(p,\delta)$-structure for some $p \in (1,2]$, and
  $\delta \in (0,\infty)$.
Let $\O\subset\setR^3$ be a bounded domain with $C^{2,1}$
  boundary and let $\ff\in L^{p'}(\O)$. Then, the unique weak solution
  $\bue\in W^{1,2}_{0}(\Omega)$ of the problem~\eqref{eq:plasticity-app}
  satisfies
  \begin{equation*}
    \begin{aligned}
      \intO   \abs{\nabla \bF(\bD\bue)}^2+\abs{\nabla \bF^{A}(\bD\bue)}^2\,d\bx
      &\le C \,, 
    \end{aligned}
  \end{equation*}  
  where $C$ depends on the characteristics of $\bS$, $\delta$, $\|\ff\|_{p'}$, $\abs{\Omega}$, and the
  $C^{2,1}$-norms of the local description of $\partial \Omega$.  In
  particular, the above estimate implies that
  $\bue$ is uniformly bounded with respect to $A\ge 1$ in $W^{2,\frac{3p}{p+1}}(\Omega).  $
\end{proposition}
\begin{proof}
  The assertion follows in the same way as in the proof of Proposition
  \ref{prop:regularity-A}. The estimate for $\nabla \bF(\bD\bue)$ in
  $L^2(\Omega)$ follows from Corollary \ref{cor:A}. It implies in turn
  that $\bue $ is bounded uniformly with respect to $A\ge 1$ in
  $W^{2,\frac{3p}{p+1}}(\Omega)$ by using \cite[Lem.~4.5]{bdr-7-5}.
\end{proof}
\section{Passing to the limit}
\label{sec:5}
The final step concerns passing to the limit $A\to\infty$. The unique solution
of~\eqref{eq:plasticity} will be obtained as
\begin{equation*}
  \bu:=\lim_{A\to\infty}\bue
\end{equation*}
with the limit taken in appropriate function spaces. 
\begin{remark}
  It will be needed to extract several sub-sequences, but we still write simply
  $A\to\infty$ to avoid using too heavy notation.
\end{remark}
By uniform --with respect to $A$-- estimates in $W^{1,2}(\Omega)$ in
Proposition \ref{prop:regularity-A-ind} it directly follows that
$\bF^{A}(\bD\bu^{A})$ has a weak limit which we denote as $\widehat{\bF}\in
W^{1,2}(\Omega)$.  Moreover, Proposition \ref{prop:regularity-A-ind}
also yields that
$\norm{\nabla^{2}\bue}_{3p/(p+1)}\leq C$, with
a constant independent of~$A$.  Hence, the compact Sobolev embedding
$W^{2,\frac{3p}{p+1}}(\Omega)\hookrightarrow\hookrightarrow
W^{1,1}(\Omega)$ implies the strong
convergence of gradients in $L^1(\Omega)$. This also implies that
$\bD\bu^{A}(\bx)\to\bD\bu(\bx)$ for almost every $\bx\in \Omega$.

Combining these two facts with $\lim_{A\to\infty}\bF^{A}(\bP)=\bF(\bP)$, which is valid
uniformly with respect to any compact set in $\setR^{3\times 3}$, and the lower semi-continuity of the norm it
follows that
\begin{gather*}
 \lim_{A\to\infty}\bF^{A}(\bD\bu^{A})=\bF(\bD\bu)\qquad\text{
      weakly in }W^{1,2}(\Omega) \text{ and a.e. in }\Omega\,,
    \\
    \leftline{\text{and also}}
    \\
    \int\limits_{\Omega}|\nabla\bF(\bD\bu)|^{2}\,d\bx\leq C\,.
\end{gather*}
Observe that $\widehat{\bF}= \bF(\bD\bu)$ since weak limit in Lebesgue spaces and the
a.e. limit coincide. 

It remains to be proved that $\bu$ is the unique solution
of~\eqref{eq:plasticity}. From the construction of $\bS^A$ it follows
$\mathbf{S}^{A}(\bP)\to\mathbf{S}(\bP)$, uniformly with respect to any
compact set in $\setR^{3\times 3}$. This fact, coupled with the almost everywhere convergence of
$\bD\bu^{A}$, implies that
\begin{equation*}
  \lim_{A\to\infty}\mathbf{S}^{A}(\bD\bu^{A}(\bx))=\mathbf{S}(\bD\bu(\bx))
  \qquad
  a.e.\ \bx\in \Omega\,, 
\end{equation*}
which is nevertheless \textit{not} enough to infer directly that
\begin{equation*}
\lim_{A\to\infty}
\int\limits_{\Omega}\mathbf{S}^{A}(\bD\bu^{A})\cdot\bD\bw\,d\bx=\int\limits_{\Omega}\mathbf{S}(\bD\bu)\cdot
\bD\bw\,d\bx\qquad \forall \,\bw\in C^{\infty}_{0}(\Omega)\,,
\end{equation*}
and to pass to the limit in the weak formulation. To this end we need
-for instance- additionally an uniform bound on
$\mathbf{S}^{A}(\bD\bu^{A})$ in $L^{q}(\Omega)$, for some $q>1$. This
would imply that
$\mathbf{S}^{A}(\bD\bu^{A})\weakto \widehat{\mathbf{S}}$ in
$L^{q}(\Omega)$, and that the limit can be identified as
$\widehat{\mathbf{S}}=\mathbf{S}(\bD\bu)$, again by the identification
of weak and almost everywhere limits in the Lebesgue spaces.

Observe that from the definition of $\bS^A$ we have
(cf.~Proposition \ref{prop:SA-ham} and  Lemma
\ref{lem:UA})  for $p \in (1, 2]$ that
\begin{equation*}
  |  \mathbf{S}^{A}(\bD\bu^{A})|\leq
 c\,\delta^{p-2}|\bD\bu^{A}|\,.
\end{equation*}
On the other hand, the estimate $ \bF(\bD\bu^{A})\in W^{1,2}(\Omega)$, which is uniform
with respect to $A$, implies by the Sobolev embedding $W^{1,2}(\Omega)\hookrightarrow
L^6(\Omega)$ that $\|\bF(\bD\bu^{A})\|_{6}\leq C $. Using the properties of $\bF$, it
follows  for $p \in (1, 2]$ that (cf.~Proposition \ref{prop:SA-ham} and Lemma~\ref{lem:UA})
\begin{equation*}
  \|\bD\bu^{A}\|_{{3p}}\leq C\,.
\end{equation*}
Hence we get that $\mathbf{S}^{A}(\bD\bu^{A})$ is bounded uniformly in 
$L^{3p}(\Omega)$ for $1<p \le 2$. This finally allows us to pass to the limit in the weak
formulation, showing that $\bu$ solves~\eqref{eq:plasticity}. Thus, we
proved Theorem \ref{thm:MT}.
\vspace*{-4mm}
\section{On the time-dependent problem}
\label{sec:parabolic}
\vspace*{-2mm}
In this section we state the natural regularity results in the
time-dependent case. These results can be proved by adapting the
method used in the steady situation. Thus, we just give the statements
of the needed results and explain necessary changes.
\begin{remark}
  Results from this section are partially contained
  in~\cite{br-parabolic}, where they are proved with a different
  approximation method and under more restrictive assumptions on the
  data, which however yield also regularity in time, i.e., 
  it is proved there that in addition 
  $\frac{\partial}{\partial t}\bF(\bD\bu) \in L^2(\Omega)$
  holds. Here, we are keeping the minimal assumptions to prove the
  natural regularity with respect to the spatial variables. However,
  additional assumptions on the data would allow to fully recover all
  the regularity results from \cite{br-parabolic}. The result
  presented in this section is the ($p\leq2$)-counterpart
  of~\cite[Thm.~3.4]{br-multiple-approx}.
\end{remark}
We now show how the initial boundary value
problem \vspace*{-2mm}
\begin{align}
  \label{eq:plasticity-parabolic}
  \begin{aligned}
\frac{\partial\bfu}{\partial t}    -\divo \bfS (\bfD\bfu) &= \bff
    \qquad&&\text{in } I\times\Omega\,,
    \\
    \bu &= \bfzero &&\text{on } I\times\partial \Omega\,,
    \\
\bu(0)&=\bu_0&&\text{in } \Omega\,,
  \end{aligned}
\end{align}
where $I=(0,T)$, for some $T>0$, can be handled by adapting the tools used in the
steady case. First we introduce the notion of a regular solution.
\begin{definition}[Regular solution]
  Let the operator $\bS=\partial\pot$
  in~\eqref{eq:plasticity-parabolic}, derived from the potential
  $\pot $, have $(p,\delta)$-structure for some $p\in(1,\infty)$, and
  $\delta\in[0, \infty)$. 
  Let $\Omega\subset\setR^3$ be a bounded domain with $C^{2,1}$
  boundary, and let $I=(0,T)$, $T\in (0,\infty)$, be a finite time
  interval. Then, we say that $\bu$ is a regular solution of
  \eqref{eq:plasticity-parabolic} if
  ${\bu \in L^p(I;W^{1,p}_0 (\Omega))}$ satisfies for all
  $\psi \in C_0^\infty (0,T)$ and all $\bw\in W^{1,p}_0(\Omega)$
\begin{equation*}
\begin{aligned}
  \int\limits _0^T\Bighskp{\frac {\partial\bu(t)}{\partial
      t}}{\bw}\,\psi(t)   + 
  \hskp{\bS(\bD\bu(t))}{\bD\bw}\,\psi (t)\,dt 
=
\int\limits
  _0^T\hskp{\bff(t)}{\bfw}\, \psi(t)\, dt\,,
\end{aligned}
\end{equation*}
and  fulfils
  \begin{align*}
    \begin{split}
      { \bfu}&\in L^{\infty}(I;W^{1,p}_{0}(\Omega))\cap {W^{1,2}(I; L^2(\Omega))}\,, 
     \\
      \bF(\bD\bfu)&\in L^{\infty}(I;L^{2}(\Omega))\cap L^{2}(I;W^{1,2}(\Omega)) \,.
    \end{split}
  \end{align*}
 \end{definition}
To formulate clearly the dependence on the data in the various
estimates  we introduce the quantity
\begin{equation*}
  |||\bu_0,\bff|||^2:=
  \int\limits_\Omega 
  |\bu_0|^2 +  |\bD\bu_0(\bx)|^p\,d\bx+\int\limits_0^T\int\limits_\Omega
  |\bff(t,\bx)|^{p'} + |\bff(t,\bx)|^{2} 
  \,d\bx \,dt\,.
\end{equation*}
As in the steady case we replace the operator $\bS$ in
\eqref{eq:plasticity-parabolic} by $\bS^A$ as in
Definition~\ref{def:SA}. The next step is the construction of a
``strong solution''. This is done by means of a Galerkin approximation
of the $A$-approximation of problem \eqref{eq:plasticity-parabolic}, 
using a priori estimates obtained by formally testing with $\bu^A$ and
$\frac{\partial \bu^A}{\partial t}$. In fact, we proceed as in the
proof of Proposition \ref{thm:existence_perturbation} and
\cite[Prop.~3.7]{br-multiple-approx}, and use also Corollary
\ref{cor:UA}. However, in contrast to the case $p >2$ in
\cite[Prop.~3.7]{br-multiple-approx} we also need to approximate the
initial condition $\bu _0$ to obtain a priori estimates independent of
$\delta^{-1}$. In fact, if we would not do so we have to handle the
term $ \function^{A}(|\bD\bu_{0}|)$ which results from testing with
$\frac{\partial \bu^A}{\partial t}$. This could be done by using
Corollary \ref{cor:UA} yielding
$\function^{A}(|\bD\bu_{0}|)\leq c\, \delta^{p-2}|\bD\bu_{0}|^{2}$,
which produces an undesired $\delta^{-1}$ dependence. We avoid this by
approximating $\bu_0$ in $W^{1,p}_0(\Omega)\cap L^2(\Omega)$ by an
$\bue_0 \in W^{1,\infty}_0(\Omega)$ satisfying
\begin{gather}
  \label{eq:app}
  \|\bD\bue_0\|_\infty \le A\,.
\end{gather}
This could be done by using the ``convolution-translation'' method,
which finds its introduction probably in the work of Puel and
Roptin~\cite{PR1967} and was re-discovered many times for different
applications to partial differential equations or simply by appealing
to standard properties of Sobolev functions and mollification. In fact,
from the proof of \cite[Thm.~5.5.2]{evans-pde} and standard
transformation and covering arguments it follows that for 
$\bu_0 \in W^{1,p}_0(\Omega)\cap L^2(\Omega)$ there exists a sequence
${(\bw_n)\subset W^{1,p}_0(\Omega)\cap L^2(\Omega)}$ and $n_0 \in
\setN$ such that
$\operatorname{supp } \bw_n \subset \Omega_{\frac 1n}:=\{\bx \in
\Omega\fdg \dist(\bx ,\partial \Omega) <\frac 1n\}$, $ n \ge n_0$,  $\|\bD\bw_n \|_{p}
\le 2\,\|\bD\bu_0\|_{p}$,  $\|\bw_n \|_{2}
\le 2\,\|\bu_0\|_{2} $, $ n \ge n_0$,  and $\bw_n \to
\bu_0$ in $W^{1,p}_0(\Omega)\cap L^2(\Omega)$. Thus, we can mollify
with a standard mollification kernel $\rho$, which yields (for
$0<\epsilon_{n}<1/2n$) a sequence 
$\bv_n:= \rho_{\vep_n} * \bw_n$ belonging  to $
C_0^{\infty}(\Omega) $ and converging to $\bu_0$  in
$W^{1,p}_0(\Omega)\cap L^2(\Omega)$. We can choose $\vep _n$ such that
it is a decreasing null sequence. Moreover, H\"older inequality yields
\begin{equation*}
  \begin{aligned}
    \| \bD \bv_n\|_{\infty}&= \| \rho_{\vep_n}*\bD \bw_n\|_{\infty}
    \\
    &\leq \smash{\frac{1}{\epsilon_n^{ 3/p}}}\|\rho\|_{p'}\|\bD\bw_{n}\|_{p}
    \le \frac{2}{\epsilon_n^{ 3/p}}\|\rho\|_{p'}\|\bD\bu_{0}\|_{p} =:A_n\nearrow \infty\,.
  \end{aligned}
\end{equation*}
For $A \in [A_n,A_{n+1})$, $ n \ge n_0$,  we set $\bue_0:=\bv_n$, which satisfies
\eqref{eq:app} and converges  to $\bu_0$  in
$W^{1,p}_0(\Omega)\cap L^2(\Omega)$.

Now we can formulate the result showing the existence of a ``strong
solution''. 
\begin{proposition}
\label{thm:existence_perturbation-parabolic}
Let the operator $\bS=\partial\pot$, derived from the potential
$\pot$, have $(p,\delta)$-structure for some $p\in(1,2]$ and
$\delta\in[0,\infty)$. Assume that
$\bu_0 \in W^{1,p}_0(\Omega)\cap L^{2}(\Omega)$ and
$\bff \in L^{p'}(I\times \Omega)$. Let $\bS^{A}$ be as in
Definition~\ref{def:SA}, and let $\bue_0$ be as constructed above,
satisfying \eqref{eq:app}. Then,  for
all $A\geq A_{n_0}$ the approximate problem
\begin{equation}
  \label{eq:eq-e-parabolic}
  \begin{aligned}
    \frac{\partial\bue}{\partial t}-\divo \bS^{A}(\bfD\bue)&=\bff
\qquad&&\text{in }I\times\Omega\,,
    \\
    \bue &= \bfzero &&\text{on } I\times\partial \Omega\,,
    \\
    \bue(0)&=\bu_{0}^{A}&&\text{in }\Omega\,,
  \end{aligned}
\end{equation}
possesses a unique strong solution
$\bue$
, i.e., $\bue\in W^{1,2}(I;L^{2}(\Omega))$ with
$\bF^{A}(\bD\bue) \in L^{\infty}(I;L^{2}(\Omega))$, which satisfies
for all
$\psi \in C_0^\infty (0,T)$ and all $\bw\in W^{1,2}_0(\Omega)$
\begin{equation*}
\begin{aligned}
  \int\limits _0^T\Bighskp{\frac {\partial\bue(t)}{\partial
      t}}{\bfw}\,\psi(t)   + 
  \hskp{\bS^{A}(\bD\bue(t))}{\bD\bfw}\,\psi (t)\,dt 
=
\int\limits
  _0^T\hskp{\bff(t)}{\bfw}\, \psi(t)\, dt\,.
\end{aligned}
\end{equation*}
In addition, the solution $\bue$ satisfies 
the estimate
\begin{gather*}
  \begin{aligned}
    &\esssup _{t \in
      I}\Big(\|\bue(t)\|_{2}^{2}+\|\bF^{A}(\bD\bue(t))\|_{2}^{2} + (\para+A)^{p-2}\norm{\nabla\bue(t)}_2^{2} +\|\bF(\bD\bue(t))\|_{2}^{2}\Big)
    \\
    &\hspace{55mm} 
    +\int\limits_{0}^{T}\Bignorm{\frac{\partial\bue(s)}{\partial
        t}}_2^2 \,ds
 \leq C\,\big (\delta^{p} +|||\bu_0,\bff|||^2\big)\,,
  \end{aligned}
 \end{gather*}
 with $C$ depending only on the characteristics of $\bS$, and $\Omega$.
\end{proposition}
\begin{proof}
  We do not give the full proof, which is a combination of 
  Proposition \ref{thm:existence_perturbation} and
 \cite[Prop.~3.7]{br-multiple-approx}. It differs from
 \cite[Prop.~3.7]{br-multiple-approx} mainly in the approximation of
 the initial condition. Thus, we just derive the a priori
 estimates.

Formally taking $\bue$ as test function
  in~\eqref{eq:plasticity-parabolic} we directly get
  \begin{equation*}
    \frac{1}{2}\|\bue(t)\|^{2}_{2}+\int\limits_{0}^{t}\int\limits_{\Omega}\function^{A}(|\bD\bue(s)|)\,d\bx
    ds\leq
    \frac{1}{2}\|\bu_{0}^{A}\|^{2}_{2}+C 
    \int\limits_{0}^{t}\int\limits_{\Omega}\function^{*}(|\bff(s)|)\,d\bx 
    ds\,,  
  \end{equation*}
where the external force is treated (for a.e. $t\in I$) as in the proof
of Proposition~\ref{thm:existence_perturbation} (note that for this estimate the
special choice of $\bue_{0}$ is not essential). 

  The second estimate is obtained by
  testing~\eqref{eq:plasticity-parabolic} by the time derivative of $\bue$. In
  this way  we get
  \begin{gather*}
    \begin{aligned}
    &\int\limits_{0}^{t}\Bignorm{\frac{\partial\bue(s)}{\partial
        t}}_2^2 \,ds+ \int\limits_{\Omega}\function^{A}(|\bD\bue(t)|)\,d\bx
    \leq c\int\limits_{\Omega}\function^{A}(|\bD\bu_{0}^{A}|)\,d\bx
    +c\int\limits_{0}^{t}\|\bff(s)\|_{2}^{2}\,ds\,.
  \end{aligned}
 \end{gather*}
The problem is that $\function^{A}$ has a quadratic growth, while
$\bD\bu_{0}$ belongs to $L^{p}(\Omega)$. To resolve this we take
advantage of the special approximation $\bue_0$. In view of
\eqref{eq:app} and the definition of $\function^A$ we get 
\begin{equation*}
  \function^{A}(|\bD\bue_{0}(\bx)|)=
  \function(|\bD\bue_{0}(\bx)|)\qquad \text{for a.e. }\bx\in \Omega\,.
\end{equation*}
Consequently, testing \eqref{eq:eq-e-parabolic} with $\bue$ and
$\frac{\partial \bue}{ \partial t}$ results in 
  \begin{gather*}
    \begin{aligned}
    &
    \frac{1}{2}\|\bue(t)\|^{2}_{2}+\int\limits_{\Omega}\function^{A}(|\bD\bue(t)|)\,d\bx+\int\limits_{0}^{t}\int\limits_{\Omega}\function^{A}(|\bD\bue(s)|)\,d\bx+\int\limits_{0}^{t}\Bignorm{\frac{\partial\bue(s)}{\partial
        t}}_2^2 \,ds
    \\
    &
\leq 
    \frac{1}{2}\|\bu_{0}^{A}\|^{2}_{2}+\int\limits_{\Omega}\function(|\bD\bu_{0}^{A}|)\,d\bx
    +C  
    \int\limits_{0}^{t}\int\limits_{\Omega}\function^{*}(|\bff(s)|)\,d\bx 
    ds+C\int\limits_{0}^{t}\|\bff(s)\|_{2}^{2}\,ds\,.
  \end{aligned}
 \end{gather*}
The assertion follows, using the estimates from
Corollary~\ref{cor:UA}, the properties of the approximation $\bue _0$,
estimate \eqref{eq:f}, and the definition of $|||\bu_{0},\bff||| $.
\end{proof}
By using the same tools employed in Section \ref{sec:3} one can prove the regularity in
the interior and for  tangential derivatives (with estimates independent of $A$). Also
regularity in normal direction follows analogously, but the estimates
depend on~$A$. 
More precisely, by adapting the translation method used in the proof of
Proposition~\ref{prop:JMAA2017-1}, the result below can be proved.
\begin{proposition}  
  \label{prop:JMAA2017-1-parabolic}
  Let the operator $\bS=\partial \pot$, derived from the potential
  $\pot$,  have $(p,\delta)$-structure for some $p\in(1,2]$,
  and $\delta\in(0,\infty)$, with characteristics $(\gamma_{3},\gamma_{4},p)$. Let
  $\Omega\subset\setR^3$ be a bounded domain with $C^{2,1}$ boundary, let
  $\bu_{0}\in W^{1,p}_{0}(\Omega)\cap L^{2}(\Omega)$ and $\bff \in L^{p'}(I\times\Omega)$.  Then, the unique
  strong solution $\bue$ of the approximate problem~\eqref{eq:eq-e-parabolic} satisfies for
  a.e. $t\in I$
  \begin{align*}
    \begin{aligned}
      \int\limits_{0}^{t}\int\limits_{\Omega} \xi_0^2
      \abs{\nabla \bF^{A}(\bD\bue)}^2
      +      \function^{A}\big( \xi_{0}^2|\nabla^{2}\bue|\big)+(\delta+A)^{p-2}
      \xi_{0}^2|\nabla^{2}\bue|^2\,d\bx \,ds
     &\le c_{0} \,,
     \\[3mm]
     \int\limits_{0}^{t} \int\limits_{\Omega} \xi^2_P \abs{\td
       \bF^{A}(\bD\bue)}^2+
     \function^{A}\big( \xi_{P}^2|\td\nabla\bue|\big)+(\delta+A)^{p-2}
     \xi_{P}^2|\td\nabla\bue|^2\,d\bx\, ds\, & \le c_{P}\,,
        \end{aligned}
  \end{align*}
  where
  $c_{0}=c_{0} (\para, |||\bu_0,\bff|||,\norm{\xi_0}_{
    1,\infty},\gamma_3,\gamma_4,p)$, while the constant related to the
  neighborhood of $P$ is such that
  $c_{P}=c_{P}
  (\para,|||\bu_0,\bff|||,\norm{\xi_P}_{1,\infty},\norm{\Grenze_{P}}_{C^{2,1}},$
  $\gamma_3,\gamma_4,p)$.
\end{proposition}
By using Proposition~\ref{prop:JMAA2017-1-parabolic} and ellipticity of $\bS^A$ we can
write, for a.e. $(t,\bx)\in I\times\Omega$, the missing partial derivatives in the normal
direction (which is locally $\be_{3}$ after a rotation of coordinates) in terms of the
tangential ones, obtaining the following result.
\begin{proposition}  
  \label{prop:JMAA2017-2-parabolic}
  Under the assumptions of Proposition \ref{prop:JMAA2017-1-parabolic} there
  exists a constant $C_1>0$ such
  that, 
  provided in the local description of the boundary there holds
  $r_P<C_1$ in $(b3)$, where $\xi_{P}(\bx)$ is a cut-off function with
  support in $\Omega_P$, then\\[-1.5mm]
  \begin{equation*}
    \begin{aligned}
\int\limits_{0}^{t}\int\limits_{\Omega} 
      \xi^2_{P} \abs{\partial _3 \bF^{A}(\bD\bue)}^2
      +\function^{A}\big(\xi^2_{P}    \abs{\partial_3 \bD\bue}\big)\,d\bx \,ds
      \le C_{A}\,,
    \end{aligned}
  \end{equation*}
  where $C_{A}=C_{A}(\delta,|||\bu_0,\bff|||,\norm{\xi_P}_{1,\infty},\norm{\Grenze_{P}}_{C^{2,1}},\gamma{}_3,\gamma{}_4,p,
  \function,A)$.
\end{proposition}

Next, we improve the estimate in the normal direction in the sense
that we will show that they are bounded uniformly with respect to the
parameter $A\geq A_{n_0}$. At this stage the time-derivative is treated as an
$L^{2}$-term on the right-hand side, while an $L^{p'}$-estimate
would be needed to estimate it properly. This can be overcome by
appropriate integration by parts.  This step involves multiplying the
equations by $\xi_{P}^{2}\partial^2_{33}\bue$ and integrating by parts
over the whole domain. To this end, the following technical result is
used to justify the treatment of the time derivative.
\begin{lemma}
  \label{lem:time-derivative}
    Let $\partial\Omega\in C^{2,1}$ and let
  $\bv\in L^{2}(I;W^{2,2}(\Omega)\cap W^{1,2}_{0}(\Omega))\cap
  W^{1,2}(I;L^{2}(\Omega))$. Then, for all $t\in[0,T]$ it holds\\[-1.5mm]
  \begin{equation*}
  -\int\limits_{0}^{t}\int\limits_{\Omega}\frac{\partial \bv}{\partial
    t}\partial_{33}^2\bv\,d\bx\,
  dt=\frac{1}{2}\|\partial_{3}\bv(t)\|_2^{2}-\frac{1}{2}\|\partial_{3}\bv(0)\|_2^{2}\,. 
\end{equation*}
\end{lemma}
Note that this result requires that $\bu_{0}\in W^{1,2}_{0}(\Omega)$
and starting from this point we need further regularity of the initial
condition.  With Lemma~\ref{lem:time-derivative} one can prove the
following result:
\begin{proposition}
  \label{prop:main-parabolic}
  Let the same hypotheses as in Proposition
  \ref{prop:JMAA2017-1-parabolic} be satisfied and assume also
  $\bu_{0}\in W^{1,2}_{0}(\Omega) $. Let the local
  description $\Grenze_{P}$ of the boundary and the localization
  function $\xi_P$ satisfy ${(b1)-(b3)}$ and $(\ell 1)$
  (cf.~Section~\ref{sec:bdr}). Then, there exists a constant $C_2>0$
  such that the unique strong solution
  $\bue $ 
  of the approximate
  problem~\eqref{eq:eq-e-parabolic} satisfies for every
  $P\in \partial \Omega$ and a.e.~$t \in I$\\[-1.5mm]
  \begin{equation*}
\int\limits_{0}^{t}    \int\limits_\Omega 
      \xi^2_P |\partial_3\bF^{A}(\bD\bue)|^2\,d\bx\, ds
     \leq C\,, 
  \end{equation*}
  provided $r_P<C_2$ in $(b3)$, with $C$ depending on the
  characteristics of $\bS$, $\delta$, 
  $|||\bu_0,\bff|||$, $\|\bD\bu_0\|_2$, $\norm{\xi_P}_{
    1,\infty},\norm{\Grenze_{P}}_{C^{2,1}}$, and $C_2$.
\end{proposition}
From Proposition \ref{thm:existence_perturbation-parabolic}
and Proposition \ref{prop:main-parabolic} we deduce in the same way as
in the proof of Proposition \ref{prop:regularity-A-ind}:

\begin{proposition}
  \label{prop:main-parabolic1}
  Under the assumption of Proposition \ref{prop:main-parabolic} the
  unique strong solution $\bue
  $ of the approximate problem~\eqref{eq:eq-e-parabolic} satisfies\\[-1mm]
  \begin{equation*}
    \esssup _{t \in
      I}\Big(\|\bue(t)\|_{2}^{2}    +\|\bF(\bD\bue(t))\|_{2}^{2}\Big) +\int\limits_{0}^{T}    \int\limits_\Omega 
    |\nabla \bF(\bD\bue)(s)|^2 +\Bigabs{\frac{\partial\bue(s)}{\partial
        t}}_2^2 \,d\bx\, ds
    \leq C\,
  \end{equation*}
  with $C$ depending on the characteristics of $\bS$,
  $\delta$, $|||\bu_0,\bff|||$, $\|\bD\bu_0\|_2$, and the
  $C^{2,1}$-norms of the local description of $\partial \Omega$. In
  particular, 
  $\bue$ is uniformly bounded with respect to $A\ge 1$ in $L^p(I;W^{2,\frac{3p}{p+1}}(\Omega)).  $
\end{proposition}

Finally, passing to the limit as $A\to\infty$ can be performed in a way similar to that
used in the steady case: Observe that the bound on the time derivative allows
us to use the Aubin-Lions lemma to infer the (space-time) convergence
\begin{equation*}
\bD\bue \to\bD\bu \qquad\text{a.e. in $ I\times \Omega$, \ \  and strongly in
  $L^2(I\times \Omega)$.}
\end{equation*}
The rest of the argument requires minor changes to prove finally the
following result:
 \begin{theorem}
  \label{thm:MT-parabolic}
  Let the operator $\bS$ in~\eqref{eq:plasticity-parabolic}, derived from a potential $\pot $,
  have $(p,\delta)$-structure for some $p\in(1,2]$, and $\delta\in(0, \infty)$.  Let $\Omega\subset\setR^3$ be a bounded domain with $C^{2,1}$
  boundary. Assume that $\bu_0\in W^{1,2}_{0}(\Omega)$ and $\bff \in L^{p'}(I\times
  \Omega)$.  Then, the system~\eqref{eq:plasticity-parabolic} has a unique regular solution
  with norms estimated only in terms of the characteristics of $\bS$,
  $\delta$,  $\Omega$,
  $\|\bu_0\|_{1,2}$, and $\|\bff\|_{p'}$.
\end{theorem}

\begin{acknowledgement}
Luigi C. Berselli was partially supported by a grant of the group
GNAMPA of INdAM.
\end{acknowledgement}
\def\cprime{$'$} \def\cprime{$'$} \def\cprime{$'$}

\end{document}